\documentclass{article}
\usepackage[utf8]{inputenc}
\usepackage[colorlinks=true,allcolors=blue]{hyperref}
\usepackage[sort&compress,numbers]{natbib}
\usepackage{graphicx}

\usepackage{amsmath,amsthm,bbm,amssymb}
\usepackage[dvipsnames]{xcolor}

\DeclareMathOperator{\im}{Im}

\DeclareMathOperator{\logg}{\log\log}

\DeclareMathOperator*{\argmax}{argmax}

\def\R{\mathbb{R}}

\def\Z{\mathbb{Z}}
\def\P{\mathbb{P}}

\def\T{\mathbb{T}}

\def\cA{\mathcal{A}}

\def\cC{\mathcal{C}}
\def\cD{\mathcal{D}}
\def\cF{\mathcal{F}}
\def\cG{\mathcal{G}}
\def\cH{\mathcal{H}}

\def\cN{\mathcal{N}}

\def\cP{\mathcal{P}}

\def\cU{\mathcal{U}}
\def\cS{\mathcal{S}}

\def\cE{\mathcal{E}}
\def\cX{\mathcal{X}}

\def\cZ{\mathcal{Z}}

\newcommand{\E}{\mathbb{E}}

\newcommand{\mean}[1] {\E\left\{{#1}\right\}}

\newcommand{\ind}{\boldsymbol{\mathbbm{1}}} 
\newcommand{\indf}[1]{\ind\set{#1}}

\newcommand{\var}[1]{\mathrm{Var}\param{{#1}}}

\newcommand{\size}[1]{\left|#1\right|}

\newcommand{\set}[1]{\left\{#1\right\}}
\newcommand{\sbrk}[1]{\left[#1\right]}

\newcommand{\param}[1]{\left(#1\right)}

\newcommand{\prob}[1]{\mathbb{P}\left(#1\right)}

\newcommand{\eps}{\epsilon}

\newcommand{\by}{\mathbf{y}}
\newcommand{\bx}{\mathbf{x}}

\newcommand{\bp}{{\bf{p}}}

\def\hr{\hat{r}}

\def\hbr{\hat{\mathbf{r}}}

\providecommand{\setthms}[1]{#1}
\setthms{
\newtheorem{lem}{Lemma}[section]
\newtheorem{thm}[lem]{Theorem}

\theoremstyle{definition}

}
\newcommand{\cech}{\v{C}ech }

\newcommand{\erdren}{Erd\H{o}s-R\'enyi }

\def\bk{\beta_{k,n}}
\def\bk{\mathbf{k}}

\newcommand{\ninf}{n\to\infty}
\newcommand{\pois}[1]{\mathrm{Pois}\param{{#1}}}

\newcommand{\limninf}{\lim_{\ninf}}

\definecolor{mygreen}{rgb}{0, 0.68, 0.31}
\definecolor{myred}{rgb}{1.0, 0,0}

\numberwithin{equation}{section}

\def\bsplit#1\esplit{\begin{split} #1 \end{split} }
\def\splitb#1\splite{\begin{split} #1 \end{split} }
\def\beq#1\eeq{\begin{equation} #1 \end{equation}}
\def\eqb#1\eqe{\begin{equation} #1 \end{equation}}

\def\hp{\hat{p}}
\def\hbp{\hat{\bp}}

\def\CM{\mathrm{CM}}
\def\SCM{\mathrm{SCM}}
\def\HSCM{\mathrm{HSCM}}

\def\mic{\mathrm{mic}}
\def\can{\mathrm{can}}
\def\hyp{\mathrm{hyp}}

\def\bla{\boldsymbol\lambda}

\def\bM{\mathbf{M}}

\def\oZ{\overline{Z}}
\def\uZ{\underline{Z}}

\title{Random Simplicial Complexes:\\Models and Phenomena}
\author{Omer Bobrowski and Dmitri Krioukov}
\date{}

\begin{document}

\maketitle

\begin{abstract}
  We review a collection of models of random simplicial complexes together with some of the most exciting phenomena related to them. We do not attempt to cover all existing models, but try to focus on those for which many important results have been recently established rigorously in mathematics, especially in the context of algebraic topology. In application to real-world systems, the reviewed models are typically used as null models, so that we take a statistical stance, emphasizing, where applicable, the entropic properties of the reviewed models. We also review a collection of phenomena and  features observed in these models, and split the presented results into two classes: phase transitions and distributional limits. We conclude with an outline of interesting future research directions.
\end{abstract}

\section{Introduction}

Simplicial complexes serve as a powerful tool in algebraic topology, a field of mathematics fathered by Poincar\'e at the end of the 19th century. This tool was built---or discovered, depending on one's philosophical view---by topologists in order to study topology. Much more recently, over the past decade or so, this tool was also discovered by network and data scientists who study complex real-world systems in which interactions are not necessarily diadic. The complexity of interactions in these systems is amplified by their stochasticity, making them difficult or impossible to predict, and suggesting that these intricate systems should be modeled by random objects. In other words, the combination of stochasticity and high-order interactions in real-world complex systems suggests that models of \emph{random} simplicial complexes may be useful models of these systems.

From the mathematical perspective, the study of random simplicial complexes combines combinatorics and probability with geometry and topology. As a consequence, the history of random simplicial complexes is quite dramatic. The drama is that this history is super-short, compared to the histories of probability and topology taken separately. For example, the first and simplest model of random simplicial complexes, the Linial-Meshulam model~\cite{linial_homological_2006}, appeared only in 2006. The main reason behind this dramatic delay is that probability and topology had been historically at nearly opposite extremes in terms of methods, skills, intuition, and esthetic (dis)likes among mathematicians. Fortunately, the wall is now dismantled, and over the last 15 years or so, the field of random topology has been growing explosively, as our review attempts to convey through the lens of random simplicial complexes.

We do not attempt to review \emph{all} existing models and results related to random simplicial complexes, which is a mission impossible. Instead, we try to focus on those models for which some exciting phenomena---which we review as well---have been rigorously established in mathematics, especially in topology. It is not entirely coincidental that a majority of these models are particulary attractive not only from the topological and probabilistic perspectives, but also from the statistical and information-theoretic perspectives. This is because these models tend to be statistically unbiased, in the sense that they are canonical models satisfying the maximum entropy principle. As a consequence, they can be used as the correct null models of real-world complexes exhibiting certain structural features.

We take this statistical stance in our review of models in Section~\ref{sec:models}. This review is then followed by the review of some of the most exciting phenomena related to these models, which we split between phase transitions (Section \ref{sec:PT}) and distributional limit theorems (Section \ref{sec:LT}). Many of the presented results are higher-dimensional analogues of the well-known phenomena in random graphs related to connectivity, giant component, cycles, etc.,
therefore we preamble, where possible, the higher-dimensional statements with brief recollections of their one-dimensional counterparts.
The focus on mathematics taken in this review precludes us unfortunately from covering many exciting subjects, such as models of growing complexes~\cite{bianconi_complex_quantum_network_manifolds_2015,bianconi_complex_quantum_network_geometries_2015,bianconi_network_2016,courtney_weighted_2017,fountoulakis_dynamical_2019,kovalenko_growing_2021}, their applications, and phenomena in them including percolation~\cite{bianconi_percolation_2019,iacopini_simplicial_2019,schaub_random_2020,gambuzza_stability_2021,lee_homological_2021,battiston_networks_2020}. Yet in the concluding Section~\ref{sec:conclusion} that outlines our view on interesting future directions, we also comment briefly on some applications and their implications for different models of random simplicial complexes.

\section{Review of models}\label{sec:models}

In application to real-world systems, the models of random simplicial complexes reviewed here are typically used as null models reproducing a particular property of interest. In other words, these models are not intended to be ``correct models'' of real-world systems, but they are intended to be \emph{correct null models} of these systems. By ``correct'' we mean here a model that is statistically unbiased, and by ``statistically unbiased'' we mean a model that maximizes entropy across all models that have a desired property. Therefore, we begin this section with a brief recollection of basic facts behind the maximum entropy principle and canonical ensembles, which we call canonical models in this review. We then observe that, with a few exceptions, the reviewed models are higher-dimensional generalizations of $1$-dimensional random simplicial complexes, which are random graphs. Therefore, we preamble, where possible, the definition of a higher-dimensional model with its $1$-dimensional counterpart. We finish the section with a short summary of the maximum entropy properties of the reviewed models.

\subsection{Maximum entropy principle and canonical models}\label{sec:maximum-entropy}

The maximum entropy principle~\cite{jaynes_information_1957} formalizes the concept of statistical unbiasedness of a null model. Indeed, Shannon entropy is the unique measure of information satisfying the basic axioms of continuity, monotonicity, and system and subset independence~\cite{shannon_mathematical_1948}, and the maximum entropy principle follows directly from these axioms coupled with the additional consistency axioms of uniqueness and representation invariance~\cite{shore_axiomatic_1980,tikochinsky_consistent_1984,skilling_axioms_1988}. For these reasons, a maximum-entropy model
is the unique model that contains all the information that the model is asked to model, and, more importantly, that does \emph{not} contain any other junk information.

Formally, let $\cX$ be a space of graphs or simplicial complexes, and $\P$ a probability distribution on $\cX$ corresponding to a \emph{model} of $\cX$: $\P(X)$ is the probability with which the model generates~$X\in\cX$. The \emph{Shannon entropy} of~$\P$ is $S(\P)=-\sum_{X\in\cX}\P(X)\log\P(X)$. Let $x_q:\cX\to\R$, $q=1,2,\ldots,Q$, be a finite collection of functions that we call \emph{properties} of~$X$, and let $y_q\in\R$ be a set of numbers that we will associate with \emph{values} of properties $x_q$. In general, properties $x_q$ can take values in spaces that are more sophisticated than~$\R$, but $\R$ suffices for us here. Denote $\bx=\set{x_q}_{q=1}^Q$, $\by=\set{y_q}_{q=1}^Q$, and let $\rho$ be a probability distribution on~$\R^Q$.

Given a pair of properties~$\bx$ and their values~$\by$, the \emph{microcanonical model} is then the one that maximizes entropy under the sharp constraints that the values of the properties $\bx$ of $X$ must be equal to $\by$ exactly: $\P_\mic(\bx,\by)=\argmax_\P\set{S(\P):\bx(X)=\by}$. This means that $\P_\mic(\bx,\by)$ is the uniform distribution over all such $X$s:
\[
\P_\mic[\bx,\by](X)=\frac{1}{\cN(\bx,\by)},
\]
where $\cN(\bx,\by)=\size{\set{X\in\cX:\bx(X)=\by}}$.

The \emph{canonical model} is the one that maximizes entropy under the soft constraints that the values of the properties~$\bx$ are equal to $\by$ in expectation: $\P_\can(\bx,\by)=\argmax_\P\set{S(\P):\E\,\bx=\by}$. If the properties~$\bx$ are sufficiently nice (e.g., satisfy certain convexity assumptions~\cite{rassoul-agha_course_2015,chatterjee_large_2017}), then the measure $\P_\can(\bx,\by)$ is a Gibbs measure
\[
\P_\can[\bx,\by](X)=\exp[-\bla(\by)\cdot\bx(X)]/Z,
\]
where $Z=\sum_{X\in\cX}\exp[-\bla(\by)\cdot\bx(X)]$ and the parameters $\bla(\by)$ solve the system of equations $\E\,\bx=-\partial\log Z/\partial\bla=\by$. The solution exists and is unique under the same niceness assumptions~\cite{rassoul-agha_course_2015,chatterjee_large_2017}. Since Gibbs measures are known as \emph{exponential families} in statistics, canonical models of random graphs are called \emph{exponential random graph models} there. Consequently, canonical models of simplicial complexes are \emph{exponential random simplicial complexes}~\cite{zuev_exponential_2015}.

Finally, the \emph{hypercanonical model} $\P_\hyp(\bx,\rho)$ is the canonical model with random $\by\sim\rho$. The measure $\P_\hyp(\bx,\rho)$ is thus the $\rho$-mixture of the Gibbs measures,
\[
\P_\hyp[\bx,\rho](X)=\int \P_\can[\bx,\by](X)\,d\rho(\by),
\]
reproducing a desired distribution~$\rho$ of the values~$\by$ of the properties~$\bx$ of~$X$.

It is important to notice that $\P_\can[\bx,\by](X)$ depends on~$X$ only via $\bx(X)$. For this reason, the properties~$\bx$ are called \emph{sufficient statistics}~\cite{cover_elements_2005}: it suffices to know $\bx(X)$ to know $\P_\can[\bx,\by](X)$; no further details about $X$ are needed. It follows that all $X$'s with the same value of $\bx(X)$ are equally likely in the canonical model, so that it is a probabilistic mixture of the corresponding microcanonical models, while the hypercanonical model is a probabilistic mixture of the canonical models.

Canonical models of random graphs and simplicial complexes tend to be more tractable than their microcanonical counterparts, explaining why the best-studied models of random simplicial complexes are based on canonical models of random graphs, versus microcanonical ones. We begin with the simplest models based on the \erdren graphs.

\subsection{Complexes based on homogeneous random graphs}\label{sec:homo-complexes}

\paragraph{The \erdren graphs $G(n,p)$ and $G(n,M)$.}
Perhaps the best-studied random graph model is the $G(n,M)$ model of random graphs with $n$ nodes and $M$ edges. All such graphs are equiprobable in the model, so that the model is microcanonical. It was introduced in 1951 by Solomonoff and Rapoport~\cite{solomonoff_connectivity_1951}, and then studied by Erd\H{o}s and R\'enyi in 1959~\cite{erdos_random_1959}.
Its canonical counterpart, introduced by Gilbert in 1959~\cite{gilbert_random_1959}, is the $G(n,p)$ ensemble of random graphs in which every possible edge between the $n \choose 2$ pairs of nodes exists independently with probability~$p$. The sufficient statistic in this canonical model is the number of edges~$M$.
If ${n\choose2}p/M\to1$ in the large graph limit $n\to\infty$, then the microcanonical $G(n,M)$ and canonical $G(n,p)$ models are asymptotically equivalent according to all definitions of such equivalence~\cite{anand_entropy_2009,squartini_breaking_2015,janson_asymptotic_2010}.

\paragraph{The Linial-Meshulam complex $Y_2(n,p)$.}
The constructive definition of the $G(n,p)$ model can be rephrased as follows: take a complete $0$-complex, which is a set of $n$ vertices, and then add $1$-simplexes (edges) to it at all possible $n\choose2$ locations independently with probability $p$. The result is a random $1$-dimensional complex~$Y_1(n,p)=G(n,p)$.

The Linial-Meshulam model~\cite{linial_homological_2006} is a straightforward $2$-dimensional analogy of $G(n,p)$: take a complete $1$-complex, which is the complete graph of size $n$, and then add $2$-simplexes (filled triangles) to it at all possible $n\choose3$ locations independently with probability $p$. The result is a random $2$-dimensional complex~$Y_2(n,p)$.

\paragraph{The $d$-dimensional Linial-Meshulam complexes $Y_d(n,p)$ and $Y_d(n,M)$.}
The Linial-Meshulam complex admits the natural generalization to any dimension $d=1,2,\ldots,n-1$~\cite{meshulam_homological_2009} in the following way. Take a complete $(d-1)$-complex, and then add $d$-simplexes to it at all possible $n\choose d+1$ locations independently with probability $p$. We denote this random complex by $Y_d(n,p)$.

The microcanonical version $Y_d(n,M)$ of canonical $Y_d(n,p)$ is also well defined. Here, we can take a complete $(d-1)$-complex,  add exactly $M$ $d$-simplexes to it chosen uniformly at random  out of all the ${{n\choose d+1}\choose M}$ possibilities. Note that $Y_1(n,M)=G(n,M)$. 
The microcanonical $Y_d(n,M)$ model gained less consideration than the canonical $Y_d(n,p)$ one, but some of its aspects were studied in~\cite{kahle_inside_2016,luczak_integral_2018}.
It should be easy to show that $Y_d(n,M)$ is asymptotically equivalent to $Y_d(n,p)$, but this has not been done.

\paragraph{The random flag complexes $X(n,p)$ and $X(n,M)$.} The Linial-Meshulam complexes are just one way to generalize the \erdren graphs to higher dimensions. Another straightforward generalization is to extend $G(n,p)$ into a random flag complex $X(n,p)$~\cite{kahle_topology_2009,kahle_sharp_2014}. Recall that the \emph{flag complex} of a graph $G$ is the simplicial complex $X$ obtained by filling all the $(k+1)$-cliques in $G$ with $k$-simplexes, for all $k=1,2,\ldots,n-1$. We can similarly define $X(n,M)$ as the flag  complex of $G(n,M)$. To the best of our knowledge, the $X(n,M)$ model has not been considered in the past.

\paragraph{The multi-parameter complexes $X(n,\bp)$ and $X(n,\bM)$.}
A further generalization of $G(n,p)$, subsuming both the Linial-Meshulam and flag complexes, is the multi-parameter complex $X(n,\bp)$
considered in~\cite{kahle_topology_2014,costa_homological_2015,costa_random_2016,costa_large_2016}. Let $\bp = (p_1,p_2,\ldots,p_{n-1})$ where $p_k \in [0,1]$, $k=1,2,\ldots,n-1$, is the probability of the existence of simplexes of dimension~$k$ in $X(n,\bp)$. Given this vector of probabilities, the complex $X(n,\bp)$ is then defined as follows. First, take $n$ vertices and add edges to them independently with probability $p_1$, i.e., generate a $G(n,p_1)$ graph. Second, considering this graph as a $1$-skeleton of a $2$-complex, go over all the triangles (i.e.~$3$-cliques) in it, and fill each of them with a $2$-simplex independently with probability $p_2$. Do not stop here, but continue in this fashion---given the $(k-1)$-skeleton, go over all the $(k+1)$-cliques in it, and fill each of them with a $k$-simplex independently with probability~$p_k$---until $k=n-1$.

Among the $G(n,p)$-based complexes discussed thus far, the $X(n,\bp)$ complex is the most general since it subsumes both the flag complex, \[X(n,p) = X(n,(p,1,1,\ldots,1)),\] and the Linial-Meshulam complex, \[Y_d(n,p) = X(n,(\underbrace{1,1,\ldots,1}_{d-1},p,\underbrace{0,0,\ldots,0}_{n-d-1})).\]

The $X(n,\bp)$ model is canonical, and to specify its sufficient statistics, we define a {\it $k$-shell} $\partial\sigma$ to be the complete $(k-1)$-dimensional boundary of a \emph{(potential)} $k$-simplex~$\sigma$ in a complex~$X$. That is, $X$ may or may not contain~$\sigma$ ($\sigma$ can be filled or empty), but if $\partial\sigma$ is a $k$-shell, then all its simplexes are filled, so that $\partial\sigma$ can be thought of as a ``{pre-$k$-simplex}'', in the sense that it is ready to be filled with~$\sigma$, but might eventually be left empty. 

As shown in~\cite{zuev_exponential_2015} (see also Section~\ref{sec:entropy}), the sufficient statistics in $X(n,\bp)$ are not only the numbers~$M_{s,k}$ of $k$-simplexes of each dimension $k$, but also the numbers~$M_{c,k}$ of $k$-shells.
Therefore, the microcanonical counterpart of $X(n,\bp)$ is $X(n,\bM)$, where $\bM=((M_{c,1},M_{s,1}),(M_{c,2},M_{s,2}),\ldots,(M_{c,n-1},M_{s,n-1}))$, $M_{c,1}={n\choose2}$, and $M_{s,1}=M$.
This model has not been considered in the past, but it is another natural higher-dimensional generalization of~$G(n,M)$, while $X(n,M)$ and $Y_d(n,M)$ are special cases of~$X(n,\bM)$.

\begin{figure}
  \centerline{\includegraphics[width=4in]{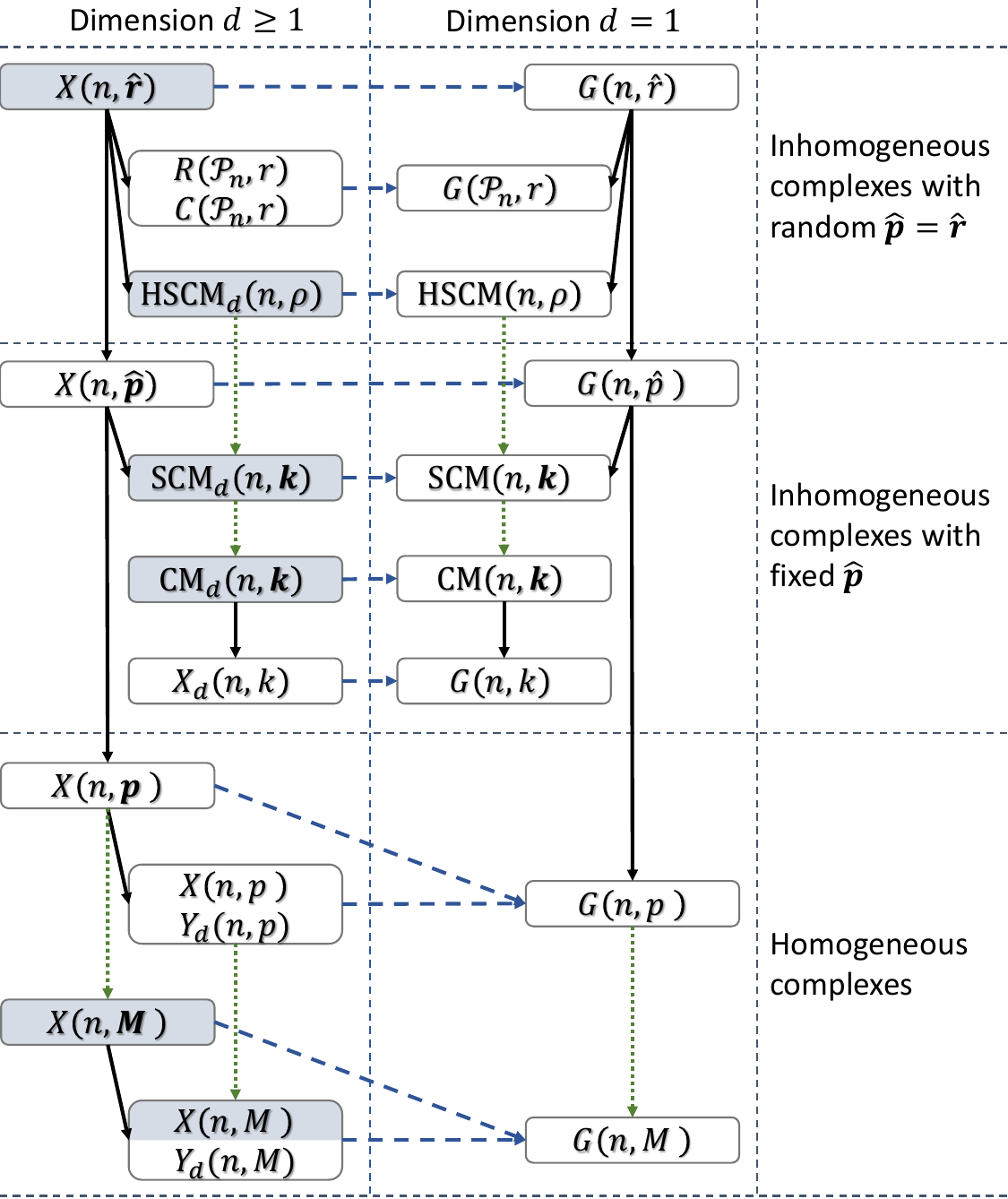}}
  \caption{{\bf The relations between the considered models of random simplicial complexes.} The vertical solid lines connect more general complexes to their special cases. The horizontal dashed lines connect higher-dimensional complexes to their $1$-dimensional cases or $1$-skeletons. The vertical dotted lines connect probabilistic mixtures (e.g., canonical models) to their constituents (e.g., microcanonical models). The shaded models have not been considered before. The review of models in Section~\ref{sec:models} proceeds roughly against the arrow directions, from the least general $G(n,M)$ to the most general $X(n,\hbr)$.}
  \label{fig:models}
\end{figure}

All other relationships between the discussed complexes are shown in Fig.~\ref{fig:models}.

\subsection{Complexes based on inhomogeneous random graphs}

\subsubsection{General complexes}

\paragraph{The inhomogeneous random graph $G(n,\hp)$.}
The edge existence probability in the random \erdren graph $G(n,p)$ does definitely not have to be the same $p$ for all edges. Each edge $i,j$ can have a different existence probability $p_{ij}$, $i,j=1,2,\ldots,n$. Such random graphs are known as inhomogeneous or generalized random graphs~\cite{soderberg_general_2002,park_statistical_2004,britton_generating_2006,bollobas_phase_2007}. We denote them by $G(n,\hp)$, where $\hp$ is the matrix of edge existence probabilities, $\hp=\{p_{ij}\}$.

\paragraph{The multi-parameter complex $X(n,\hbp)$.}
A straightforward generalization of $G(n,\hp)$ to higher dimensions was introduced in~\cite{zuev_exponential_2015}. It is also a generalization of $X(n,\bp)$, and it is defined as follows. Let
$\hbp=(\hp_1,\hp_2,\ldots,\hp_{n-1})$
be a vector of simplex existence probabilities for all possible simplexes of all possible dimensions $k=1,2,\ldots,n-1$: $\hp_1=\{p_{ij}\}$, $\hp_2=\{p_{ijl}\}$, and $\hp_k=\{p_{\sigma_k}\}$ collects the $n\choose k+1$ existence probabilities for all possible $n\choose k+1$ simplexes $\sigma_k$ of dimension $k$.
Given such a vector of tensors $\hbp$, the random complex $X(n,\hbp)$ is then generated similarly to $X(n,\bp)$: starting with an inhomogeneous random graph $G(n,\hp_1)$, go over all the $3$-cliques ($2$-shells) in this graph, and fill them with $2$-simplexes $\sigma_2$ independently with probability~$p_{\sigma_2}$. Proceed to higher dimensions $k=3,4,\ldots,n-1$ in a similar fashion, filling $k$-shells with $k$-simplexes independently with probabilities~$p_{\sigma_k}$.
The complex $X(n,\bp)$ is a special case of $X(n,\hbp)$---the one with $\hp_k\equiv p_k$.

\subsubsection{Configuration models}\label{sec:scm}

\paragraph{The configuration models $\mathrm{SCM}(n,\bk)$ and $\mathrm{CM}(n,\bk)$.}
The inhomogeneous random graphs $G(n,\hp)$ are quite general, subsuming many other important models of random graphs. In particular, they encompass the soft configuration model $\SCM(n,\bk)$~\cite{soderberg_general_2002,park_statistical_2004,britton_generating_2006,bollobas_phase_2007,bianconi_entropy_2008,garlaschelli_maximum_2008,chatterjee_random_2011}, which is the canonical model of random graphs with a given sequence of expected degrees $\bk=(k_1,k_2,\ldots,k_n)$, $k_i\geq0$, $i=1,2,\ldots,n$. The $\SCM(n,\bk)$ model is the $G(n,\hp)$ with $\hp$ given by
\beq\label{eq:scm-p}
p_{ij}=\frac{1}{e^{\lambda_i+\lambda_j}+1},
\eeq
where the parameters $\lambda_i$'s (known as Lagrange multipliers) solve the system of $n$ equations given by
\beq\label{eq:scm-k}
\sum_{j}p_{ij}=k_i.
\eeq
This equation guarantees the expected degree of node~$i$ to be $k_i$.

The sufficient statistics in $\SCM(n,\bk)$ are the degrees of all nodes in a graph. Therefore, $\SCM(n,\bk)$'s microcanonical counterpart, the
configuration model $\CM(n,\bk)$~\cite{bender_asymptotic_1978,molloy_critical_1995}, is the uniform distribution over all the graphs with the degree sequence $\bk$.
Note that $\bk$ can be any sequence of nonnegative real numbers in $\SCM(n,\bk)$, but in $\CM(n,\bk)$, $\bk$ is a graphical sequence of nonnegative integers. A sequence $\bk$ is called \emph{graphical}, if there exists a graph whose degree sequence is~$\bk$. The necessary and sufficient conditions for~$\bk$ to be graphical are the Erd\H{o}s-Gallai conditions~\cite{erdHos1960evolution}.

An important special case of $\CM(n,\bk)$ is the random $k$-regular graph $G(n,k)=\CM(n,(k,k,\ldots,k))$.

\paragraph{The $d$-dimensional configuration models $\mathrm{SCM}_d(n,\bk)$ and $\mathrm{CM}_d(n,\bk)$.}
Recall that the degree $k_i$
of a vertex ($0$-simplex)~$i$ is the number of edges ($1$-simplexes) that contain~$i$. In a similar vein, the {\it degree} $k_\sigma$ of $d$-simplex $\sigma$ is defined to be the number of $(d+1)$-simplices that contain~$\sigma$. 

The $d$-dimensional soft configuration model $\SCM_d(n,\bk)$ is a generalization of $\SCM(n,\bk)$, where now $\bk=\{k_\tau\}$ is a sequence of $n\choose d$ expected degrees $k_\tau\geq0$ of all  $(d-1)$-simplexes $\tau$. To construct $\SCM_d(n,\bk)$ we take the complete $(d-1)$-skeleton on $n$ vertices, and  add all possible $d$-simplexes $\sigma$  independently with probability
\beq\label{eq:scmd-p}
p_\sigma=\frac{1}{e^{\sum_\tau\lambda_\tau\indf{\tau<\sigma}}+1},
\eeq
where the summation is  over all  $(d-1)$-faces $\tau$ of $\sigma$ (as implied by the standard $\tau<\sigma$ notation), and the parameters $\lambda_\tau$ solve the system of  $n\choose d$ equations
\beq\label{eq:scmd-k}
\sum_{\sigma}p_\sigma\indf{\tau<\sigma}=k_\tau.
\eeq
Note that $\SCM_d(n,\bk)$  is a special case of $X(n,\hbp)$:
\[\SCM_d(n,\bk)=X(n,(\underbrace{1,1,\ldots,1}_{d-1},\hp_d,\underbrace{0,0,\ldots,0}_{n-d-1})),\] where $\hp_d=\{p_\sigma\}$ is given by~\eqref{eq:scmd-p}. If $d=1$, then Eqs.~(\ref{eq:scmd-p},\ref{eq:scmd-k}) reduce to Eqs.~(\ref{eq:scm-p},\ref{eq:scm-k}), respectively, and we have $\SCM_1(n,\bk)\equiv\SCM(n,\bk)$.

The $\SCM_d(n,\bk)$ is the canonical model of random complexes whose sufficient statistics are the degrees of all the $(d-1)$-simplexes in a complex. Therefore, $\SCM_d(n,\bk)$'s microcanonical counterpart, the configuration model $\CM_d(n,\bk)$, is the uniform distribution over complexes with the complete $(d-1)$-skeleton and the degree sequence of $(d-1)$-simplexes equal to~$\bk$.
In $\CM_d(n,\bk)$, $\bk$ is a realizable sequence of $n\choose d$ nonnegative integers, versus any sequence of $n\choose d$ nonnegative real numbers in $\SCM_d(n,\bk)$. The conditions for~$\bk$ to be realizable, analogous to the Erd\H{o}s-Gallai graphicality conditions in the $d=1$ case, are at present unknown.

An important special case of the configuration model is the random $k$-regular complex $X_d(n,k)=\CM_d(n,(k,k,\ldots,k))$ in which all $(d-1)$-simplexes have degree $k$ \cite{lubotzky_random_2019}. For $k=1$, such a complex is also known as an $(n,d)$-Steiner system, whose randomized construction is due to \cite{keevash2014existence}.

With the exception of~$X_d(n,k)$, the $d$-dimensional configuration models $\mathrm{CM}_d(n,\bk)$ and $\mathrm{SCM}_d(n,\bk)$ have not been considered in the past. The configuration models defined and studied in~\cite{courtney_generalized_2016,young_construction_2017} are very different and unrelated to any other model considered above. We review them in Section~\ref{sec:zcm}.

\subsection{Complexes based on inhomogeneous random graphs with random connection probabilities}

\subsubsection{General complexes}\label{sec:general-random}

\paragraph{The inhomogeneous random graph with random connection probabilities $G(n,\hr)$.}
Note that the edge probabilities $p_{ij}$ in inhomogeneous random graphs $G(n,\hp)$ can  themselves  be random~\cite{caldarelli_scale_2002,boguna_class_2003,britton_generating_2006,bollobas_phase_2007}. In that case, we replace $\hp$ with $\hr=\{r_{ij}\}$ which is a random matrix with entries in $[0,1]$, and the corresponding inhomogeneous random graph is denoted $G(n,\hr)$. This class of random graphs is extremely general and subsumes a great variety of many important random graph models, for example $G(n,\hp)$.

Of a particular interest is the case where \[r_{ij}=p(X_i,X_j),\] 
where $\{X_1,\ldots,X_n\}$ is a set of random variables in some measure space $S$, and  $p:S\times S\to [0,1]$ is a fixed function called the {\it connection probability function}. This class of models is quite general and  subsumes, for instance, the stochastic block models~\cite{holland_stochastic_1983}, hidden variable models~\cite{caldarelli_scale_2002,boguna_class_2003}, latent space models~\cite{sorokin_social_1927,mcfarland_social_1973,hoff_latent_2002}, random geometric graphs~\cite{gilbert_random_1961,penrose_random_2003}, and many other popular models. If the $X_i$'s are independent random variables uniformly distributed in $[0,1]$, and $p$ is a symmetric integrable function, then $p$ is also known as a {\it graphon} in the theory of graph limits~\cite{lovasz_large_2012,janson_graphons_2013}.

\paragraph{The multi-parameter complex $X(n,\hbr)$.}
The most general class of models considered in this chapter
is a version of the multi-parameter complex $X(n,\hbp)$ with random simplex existence probabilities~$\hbp$. To emphasize their randomness, we denote this complex by $X(n,\hbr)$, where $\hbr$ is random~$\hbp$. For maximum generality, the probabilities of existence of different simplexes, including simplexes of different dimensions, do not have to be independent random variables, so that the space of $X(n,\hbr)$ is parameterized by joint probability distributions over the vector of tensors~$\hbr$, equivalent to all nonempty subsets of~$\set{1,\ldots, n}$. As evident from Fig.~\ref{fig:models}, all other models of random complexes and graphs are special cases of $X(n,\hbr)$. In particular, $X(n,\hbp)$ is a special degenerate case of $X(n,\hbr)$. At this level of generality, the complex $X(n,\hbr)$ has not been considered or defined before, yet it is simply a higher-dimensional version of the well-studied~$G(n,\hr)$.

\subsubsection{Hypersoft configuration model}

\paragraph{The hypersoft configuration model $\HSCM(n,\rho)$.}
The hypersoft configuration model is the hypercanonical model of random graphs with a given expected degree distribution~$\rho$~\cite{caldarelli_scale_2002,boguna_class_2003,britton_generating_2006,bollobas_phase_2007,chatterjee_random_2011,hoorn_sparse_2018}. It is defined as the $\SCM(n,\bk)$ in which the expected degree sequence~$\bk$ is not fixed but random: the expected degree $k_i$ of every vertex $i$ is an independent random variable with distribution~$\rho$, $k_i\sim\rho$. In other words, to generate an $\HSCM(n,\rho)$ graph, one first samples $k_i$s independently from the distribution~$\rho$, then solves the system of equations~\eqref{eq:scm-k} to find all $\lambda_i$s, and finally creates edges with probabilities~$p_{ij}$~\eqref{eq:scm-p}. The $\HSCM(n,\rho)$ is thus a probabilistic mixture of $\SCM(n,\bk)$'s with random $\bk\sim\rho^n$. The degree distribution in $\HSCM(n,\rho)$ is the mixed Poisson distribution with mixing~$\rho$~\cite{britton_generating_2006,bollobas_phase_2007}.

An equivalent definition of $\HSCM(n,\rho)$ is based on the observation that the distribution~$\rho$ defines the joint distribution~$\Psi$ of $\Lambda=\{\lambda_{i}\}$ via~\eqref{eq:scm-k}. This means that an equivalent procedure to generate an $\HSCM(n,\rho)$ graph is to sample $\Lambda$ directly from $\Psi$ first, and then create edges with probabilities~$p_{ij}$~\eqref{eq:scm-p}. This definition makes it explicit that the $\HSCM(n,\rho)$ is a special case of $G(n,\hr)$.

For a given $\rho$, it may difficult to find the explicit form of the distribution~$\Psi$. However, in many cases with important applications---including the Pareto $\rho$ with a finite mean, for instance---it has been shown that the canonical connection probability~\eqref{eq:scm-p} and its classical limit approximation
\beq\label{eq:scm-pCL}
p_{ij}^{\mathrm{CL}}=\min\left(1,e^{-\lambda_i}e^{-\lambda_i}\right)=\min\left(1,\frac{k_ik_j}{\bar{k}n}\right),
\eeq
where $\bar{k}$ is the mean of $\rho$, lead to asymptotically equivalent random graphs~\cite{janson_asymptotic_2010}. In such cases, the random Lagrange multipliers $\lambda_i$'s are asymptotically independent, $\Psi=\psi^n$, $\lambda_i\sim\psi$, with the distribution~$\psi$ defined by the distribution~$\rho$ via the change of variables $\lambda_i =\ln\left(\sqrt{\bar{k}n}/k_i\right)$ where $k_i\sim\rho$. Generating such a graph is extremely simple: first sample either the $k_i$'s or $\lambda_i$'s independently from the distribution~$\rho$ or $\psi$, respectively, and then generate edges with probabilities~$p_{ij}^{\mathrm{CL}}$~\eqref{eq:scm-pCL}.

\paragraph{The $d$-dimensional hypersoft configuration model $\HSCM_d(n,\rho)$.}
A straightforward generalization of $\HSCM(n,\rho)$ to dimension $d$ is achieved by defining the $\HSCM_d(n,\rho)$ as the probabilistic mixture of $\SCM_d(n,\bk)$ with random $\bk\sim\rho^{n\choose d}$: the expected degree $k_\tau$ of every $(d-1)$-simplex $\tau$ is an independent random variable with distribution~$\rho$. In other words, to generate a random $\HSCM_d(n,\rho)$ complex, one first prepares the complete $(d-1)$-complex of size $n$, then samples the $n\choose d$ expected degrees~$k_\tau$ of all the $(d-1)$-simplexes~$\tau$ independently from the distribution~$\rho$, then solves the system of equations~\eqref{eq:scmd-k} to find all the $\lambda_\tau$'s, and finally creates the $d$-simplexes $\sigma$ independently with probability~$p_\sigma$~\eqref{eq:scmd-p}. For the same reasons as in the $d=1$ case, the $\HSCM_d(n,\rho)$ model is a special case of $X(n,\hbr)$ for any $d\geq1$.
The model has not been previously considered, except the $d=1$ case.

\subsubsection{Geometric complexes}

An important class of random complexes are geometric complexes, which are a special case of $X(n,\hbr)$ whose $1$-skeletons are random geometric graphs. Their randomness comes from random locations of vertices in a space, and this randomness is often modeled by either binomial or Poisson point processes whose definitions we recall next.

Let $S$ be a metric space with a probability distribution $P$ on it. The definitions below apply to any sufficiently nice pair of $S$ and $P$. However, for the sake of simplicity, in the subsequent sections of this chapter, the space $S$ will always be a flat $d$-dimensional torus $S = \T^d = [0,1]^d / \set{0\sim 1}$, while $P$ will be its Lebesgue measure~$\lambda$, i.e., $P(A)=\lambda(A)/\lambda(\T^d)=\lambda(A)$, where $\lambda(A)$ is the Euclidean volume of~$A\subseteq\T^d$.

\paragraph{The binomial point process} is simply a set of $n$ random points $\cX_n = \set{X_1,\ldots,X_n} \subset S$ sampled independently from $P$. In the simplest case $S=\T^1$, we simply sample $X_i$'s from the uniform distribution, $X_i \sim \cU(0,1)$, viewed as a circle. The distribution of the number of points in any measurable subset $A\subseteq S$ is binomial with mean~$nP(A)$, giving the process its name.

\paragraph{The Poisson point process} is the binomial process with a random number of points $N$ sampled from the Poisson distribution with mean $n$, i.e.~$\cP_n = \cX_N$, where $N\sim \pois{n}$. The distribution of the number of points in any measurable $A\subseteq S$ is Poisson with mean $nP(A)$, where $n$ is the {\it rate} of the process.
Another key property of the Poisson process is \emph{spatial independence}: the numbers of points in any pair of  disjoint subsets of $S$ are independent random variables. This property, absent in the binomial process, makes the Poisson process favorable in probabilistic analysis.

In the $n\to\infty$ limit, the two process are asymptotically identical---the binomial process converges to the Poisson process in a proper sense~\cite{last_lectures_2017}. For this reason, and given that the Poisson process is easier to deal with, we will restrict ourselves in this chapter to the Poisson process $\cP_n$ on the torus $\T^d$.
 These setting turn out to provide the most elegant results, while not losing much in terms of generality.

\paragraph{The random geometric graph~$G(\cP_n,r)$.}
We start with $G(\cX_n,r)$, which is 
 a special case of $G(n,\hr)$,
in which
\[
r_{ij}=p(X_i,X_j)=\indf{d(X_i,X_j)\le r},
\]
where $d(X_i,X_j)$ is the distance between $X_i,X_j\in\cX_n$ in $\T^d$, and $r>0$ is the connectivity radius parameter~\cite{gilbert_random_1961,penrose_random_2003}. In other words, the vertices of this random graph are a realization of the binomial process, and two vertices are connected if the distance between them in $\T^d$ is at most $r$. 
To generate the random geometric graph $G(\cP_n,r)$ we first sample $N\sim\pois{n}$, and then generate $G(\cX_N,r)$ as described above.

\paragraph{The random Vietoris-Rips complex~$R(\cP_n,r)$.}
This is the flag complex over the random geometric graph $G(\cP_n,r)$. It is thus a straightforward higher-dimensional generalization of $G(\cP_n,r)$.
Note that we can consider $R(\cX_n,r)$ for the binomial process as well, which is a special case of~$X(n,\hbr)$. Then $R(\cP_n,r) = R(\cX_N,r)$ with $N\sim\pois{n}$.

\paragraph{The random \cech complex~$C(\cP_n,r)$.}
This is a different higher-dimensional generalization of $G(\cP_n,r)$. The $C(\cP_n,r)$ rule is: draw balls $B_{r/2}(X_i)$ of radius $r/2$ around each of the points $X_i\in\cP_n$, then look for all the intersections of these balls, and for every $(k+1)$-fold intersection of the balls, add the corresponding $k$-simplex to the complex. Note that the binomial version $C(\cX_n,r)$ can still be considered as a special case of $X(n,\hbr)$ with a different rule for the creation of higher-dimensional simplexes, compared to~$R(\cX_n,r)$.

\begin{figure}
\centering
\includegraphics[width=0.4\textwidth]{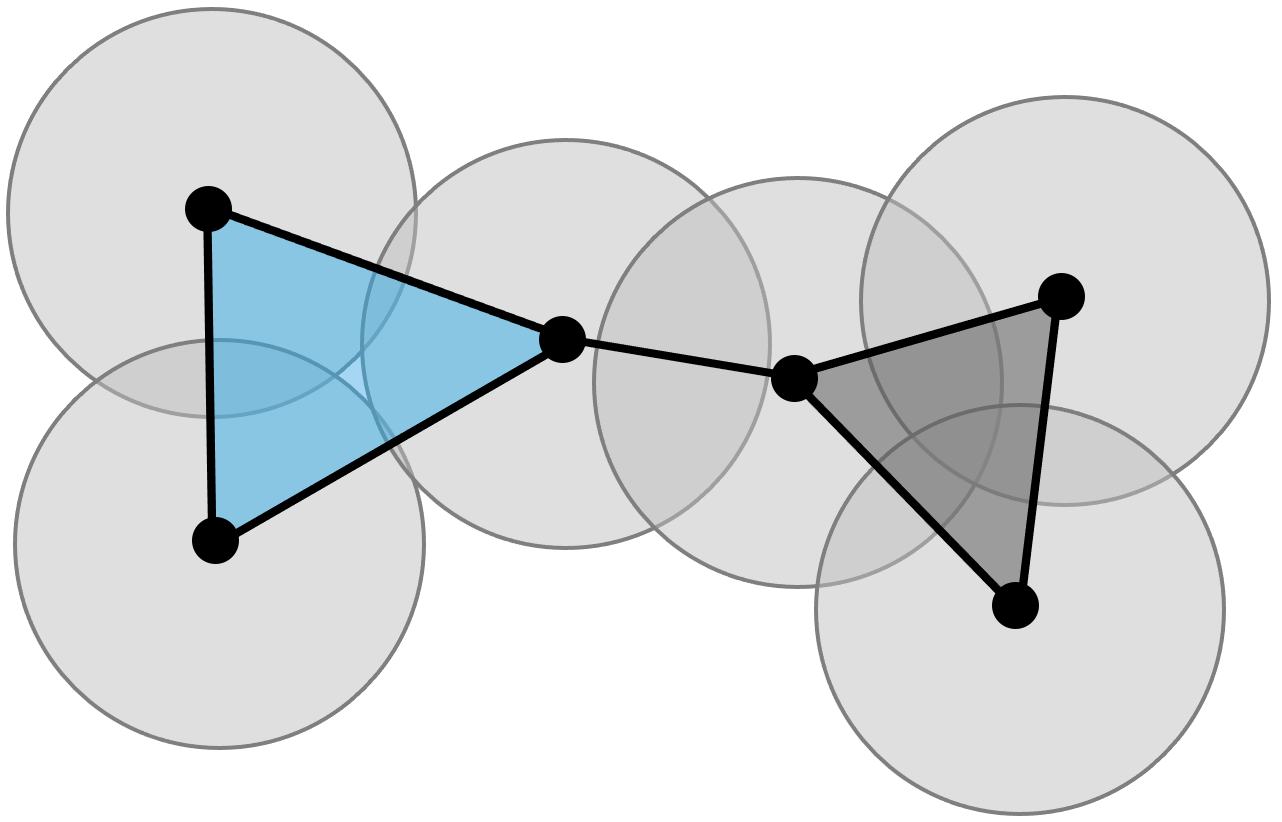}
\caption{The \cech and Rips complexes. We draw balls of a fixed radius $r$ around the points, and consider their intersection. The edges in both cases correspond to the geometric graph. While both triangles belong to the Rips complex,  the left triangle is not included in the \cech complex. \label{fig:cech_rips}}
\end{figure}

While the $1$-skeleton of both $R(\cP_n,r)$ and $C(\cP_n,r)$ is the same random geometric graph $G(\cP_n,r)$, higher-dimensional simplexes in these complexes are in general different, as illustrated in Fig.~\ref{fig:cech_rips}.
The following useful relation was proved in \cite{de_silva_coverage_2007} for any $\alpha \le \sqrt{\frac{d+1}{2d}}$:
\[
    R(\cP_n, \alpha r) \subset C(\cP_n,r) \subset R(\cP_n,r).
\]

A powerful property of the \cech complex is due to the Nerve Lemma \cite{borsuk_imbedding_1948} stating that under mild conditions we have the following isomorphism between the homology groups~$H_k$:
\[
    H_k(C(\cP_n,r)) \cong H_k(B(\cP_n,r)),
\]
where $B(\cP_n,r) = \bigcup_{X_i\in \cP_n} B_{r/2}(X_i)$.
This property, absent in the Rips complex, is very useful in the probabilistic analysis of the \cech complex, allowing one to switch back-and-forth between combinatorics and stochastic geometry.

\subsection{Complexes based on random hypergraphs}\label{sec:z-complexes}

\subsubsection{General complexes}\label{sec:lower-upper}

The main reason why the Linial-Meshulam model $Y_d(n,p)$ is defined on top of the complete $(d-1)$-skeleton, is simplicity. The complete $(d-1)$-complex underlying the $d$-complex $Y_d(n,p)$ simplifies its topological analysis drastically. Much less simple, but also much more general, is the $X(n,\hbp)$ complex in which the probability of existence of simplexes of different dimensions are all different, but this complex is still relatively simple from the probabilistic perspective, since simplexes of higher dimensions are created in a conditionally independent manner, with the conditions being the existence of required simplexes of lower dimensions. Very recently, this simplicity was sacrificed even further for the sake of even greater generality in a class of models based on random hypergraphs~\cite{farber_random_2019,farber_random_2020}, which we briefly review next.

\paragraph{The inhomogeneous random hypergraph $H(n,\hbp)$.} This is a straightforward generalization of the inhomogeneous random graph $G(n,\hp)$ to hypergraphs: every possible hyperedge~$\sigma$ (which is a nonempty subset of $\set{1,\ldots, n}$) in $H\sim H(n,\hbp)$ exists independently with probability~$p_\sigma\in\hbp$. Note that $H$ is not a simplicial complex with high probability, unless $\hbp$ is specially designed for $H$ to be a complex.

\paragraph{The lower and upper complexes $\uZ(n,\hbp)$ and $\oZ(n,\hbp)$.} Introduced in \cite{cooley2020vanishing, farber_random_2019}, these are both based on the~$H(n,\hbp)$. The lower complex $\uZ\sim \uZ(n,\hbp)$ is the largest simplicial complex that the random $H\sim H(n,\hbp)$ contains, while the upper complex $\oZ\sim \oZ(n,\hbp)$ is the smallest simplicial complex that contains the random $H\sim H(n,\hbp)$, so that $\uZ \leq H \leq \oZ$.
In other words, a simplex $\sigma$ is included in the lower complex $\uZ$ if and only if $\sigma \in H$ and also  all its faces are  hyperedges of~$H$. The other way around, every hyperedge $\sigma\in H$ is included in the upper complex $\oZ$ as a simplex together with all its faces, even if some of these faces do not happen to be  in $H$.

Compared to the graph-based models reviewed in the previous sections, the hypergraph-based models
are much more difficult to analyze, primarily because the distributions of the skeletons are highly nontrivial and carry intricate dependency structure. However, an exciting fact about these two models is that they are dual in some strict topological sense (Alexander duality), and therefore statements about one of the models can be translated with ease to  corresponding statements about the other~\cite{farber_random_2019,farber_random_2020}.

At this point, it may be tempting to proceed similarity to the previous $X$-sections, and let the hyperedge existence probabilities $\hbp$ be random $\hbr$, leading to $H(n,\hbr)$, $\uZ(n,\hbr)$, and $\oZ(n,\hbr)$. However, nothing new is gained by doing so, since the $X(n,\hbr)$ model is already general enough as it includes any possible joint distribution of simplex existence probabilities, including the distributions describing $\uZ(n,\hbp)$ and $\oZ(n,\hbp)$. That is, the lower and upper models $\uZ(n,\hbp)$ and $\oZ(n,\hbp)$ with nonrandom $\hbp$ are special cases of the more general $X(n,\hbr)$ model with the joint distribution of $\hbr$ set equal to the joint distribution of simplex existence probabilities in $\uZ(n,\hbp)$ and $\oZ(n,\hbp)$, while the lower and upper models $\uZ(n,\hbr)$ and $\oZ(n,\hbr)$ with random $\hbr$ are equivalent to $X(n,\hbr_*)$ with matching joint distributions of $\hbr_*$.

\subsubsection{$Z$-configuration models}\label{sec:zcm}

One interesting special case of the $Z$-complexes from the previous section is the $Z$-version of the $d$-dimensional (soft) configuration model considered in~\cite{courtney_generalized_2016}.

\paragraph{The $Z$-configuration models $\mathrm{Z}$-$\mathrm{CM}_d(n,\bk)$ and $\mathrm{Z}$-$\mathrm{SCM}_d(n,\bk)$.}
The {\it $d$-degree} $k_\sigma(d)$ of a $d'$-simplex~$\sigma$ is defined to be the number of simplexes of dimension $d>d'$ that contain~$\sigma$~\cite{bianconi_complex_quantum_network_manifolds_2015}, and the $Z$-configuration models are defined by a vector $\bk=\{k_i(d)\}$ of a given sequence of (expected) $d$-degrees of vertices $i=1,2,\ldots,n$, as opposed to the vector of (expected) $d$-degrees of $(d-1)$-simplexes in the $\mathrm{(S)CM}_d(n,\bk)$ in Section~\ref{sec:scm}. No low-dimensional skeleton is formed in the Z-$\SCM_d(n,\bk)$, while the probability of $d$-simplex $\sigma$ is
\beq\label{eq:scmg-p}
p_\sigma=\frac{1}{e^{\sum_i\lambda_i\indf{i<\sigma}}+1},
\eeq
where the $n$ parameters $\lambda_i$ of vertices $i$ solve the system of $n$ equations
\beq\label{eq:scmg-k}
\sum_{\sigma}p_\sigma\indf{i<\sigma}=k_i(d).
\eeq
If simplex $\sigma$ is added to the complex, then so are all its lower-dimensional faces. The Z-$\SCM_d(n,\bk)$ is a special case of the $\bar Z(n,\hbp)$ with $\hbp=(0,\ldots,0,\hp_d,0,\ldots,0)$ and $\hp_d=\{p_\sigma\}$ given by~\eqref{eq:scmg-p}.

The model is canonical with the $d$-degree of all vertices being its sufficient statistics, so that the corresponding microcanonical model Z-$\CM_d(n,\bk)$ is the uniform distribution over all the $d$-complexes with the $d$-degree sequence of the vertices equal to~$\bk$. Note that these complexes satisfy the additional constraint: they do not contain any simplex of dimension $d'<d$ that is not a face of a $d$-simplex. The hypercanonical model Z-$\HSCM_d(n,\rho)$ is the Z-$\SCM_d(n,\bk)$ with random $k_i(d)\sim\rho$.

\subsection{Random simplicial complexes as canonical models}\label{sec:entropy}

Table~\ref{tab:canonical} summarizes what is known concerning the reviewed complexes as canonical models. To convey this knowledge succinctly, it suffices to focus only on the canonical models {\it per se}, since their microcanonical constituents and hypercanonical mixtures are coupled to them as discussed in Section~\ref{sec:maximum-entropy}.

\begin{table}
  \centering
  \begin{tabular}{|l|l|l|l|l|}
    \hline
    Models & Spaces $\cS$ & Properties $\bx$ & Refs & Comments  \\ \hline
    $G(n,p)$, $G(n,M)$ & $\cG_n$ & $M$ & \cite{park_statistical_2004} & \eqref{com:match} \\
    $X(n,p)$, $X(n,M)$ & $\cF_n$ & $M$ & \cite{zuev_exponential_2015} & \eqref{com:match} \\
    $Y_d(n,p)$, $Y_d(n,M)$ & $\cC_{n,d}$ & $M_{s,d}$ & \cite{zuev_exponential_2015} & \eqref{com:match} \\
    $X(n,\bp)$, $X(n,\bM)$ & $\cC_n$ & $\bM=\set{M_{c,d},M_{s,d}}_{d=1}^{n-1}$ & \cite{zuev_exponential_2015} & (\ref{com:match},\ref{com:surprise}) \\
    $G(n,\hp)$, $G(n,\hr)$ & $\cG_n$ & $\hat{\sigma}_1=\set{\sigma_1}$ & \cite{park_statistical_2004} & \eqref{com:match} \\
    $X(n,\hbp)$, $X(n,\hbr)$ & $\cC_n$ & $\hat{\boldsymbol\sigma}=\set{\set{\partial\sigma_{d}},\set{\sigma_{d}}}_{d=1}^{n-1}$ & \cite{zuev_exponential_2015} & (\ref{com:match},\ref{com:surprise}) \\
    $\mathrm{((H)S)CM}(n,\cdot)$ & $\cG_n$ & $\bk = \set{k_i}$ & \cite{park_statistical_2004} & \eqref{com:hscm} \\
    $\mathrm{((H)S)CM}_d(n,\cdot)$ & $\cC_{n,d}$ & $\bk = \set{k_{\sigma_{d-1}}}$ & here & \eqref{com:hscm_d} \\
    Z-$\mathrm{((H)S)CM}_d(n,\cdot)$ & $\cZ_{n,d}$ & $\bk = \set{k_i(d)}$ & \cite{courtney_generalized_2016} & \eqref{com:z-hscm_d} \\
    \hline
  \end{tabular}
  \caption{{\bf Canonical models of random simplicial complexes.} The first column lists the canonical models from Section~\ref{sec:models}, along with their micro- and hyper-canonical counterparts (if any), whose entropy maximization properties have been established. The second and third columns document the constraints (space of complexes and their properties) under which the model entropy is maximized. The fourth column contains references to where this maximization has been established and to further details. The last column refers to pertinent comments in the text.}\label{tab:canonical}
\end{table}

To define a canonical model, one needs to specify not only the sufficient statistics~$\bx$ (whose expectations can take any admissible values $\by$ in a particular canonical model), but also the space $\cS$ of allowed complexes over which the entropy of a canonical model is maximized. Table~\ref{tab:canonical} uses the following notations for such spaces $\cS$ of \emph{labeled} graphs and complexes:
\begin{itemize}
  \item $\cG_n$: all graphs of size~$n$;
  \item $\cF_n$: flag complexes of size~$n$;
  \item $\cC_n$: all complexes of size~$n$;
  \item $\cC_{n,d}$: $d$-complexes of size~$n$ with a complete $(d-1)$-skeleton;
  \item $\cZ_{n,d}$: $d$-complexes of size~$n$ whose all $d'$-simplexes, for all $d'<d$, are faces of $d$-simplexes.
\end{itemize}
For the sufficient statistics, Table~\ref{tab:canonical} uses the following notations that rely on the definition of a $d$-shell, which can be found at the end of Section~\ref{sec:homo-complexes}:
\begin{itemize}
  \item $M_{c,d}$ and $M_{s,d}$: numbers of $d$-shells and $d$-simplexes ($M_{s,1}=M$, the number of edges);
  \item $\partial\sigma_d$ and $\sigma_d$: $d$-shell and $d$-simplex ($\sigma_0=i$, a vertex; $\sigma_1=\set{i,j}$, an edge);
  \item $k_{\sigma_d}$ and $k_{\sigma_d}(d')$: $(d+1)$-degree and generalized $d'$-degree of $\sigma_d$ ($d<d'<n$).
\end{itemize}
Finally, the comments referred to in Table~\ref{tab:canonical} are as follows:
\begin{enumerate}
  \item\label{com:match} The values~$\by$ of properties~$\bx$ are not shown, but they are assumed to match properly the model parameters. That is, in a general canonical model we have that $\E\,\bx=\by$, so that the $G(n,p)$ model, for example, is the canonical model with $\bx = M$ and $\by = \E\,M = {n\choose2}p$. Similarly, the $X(n,\hbp)$ model is the canonical model with $\E\,\sigma_d=p_{\sigma_d}\E\,\partial\sigma_d$ and $\E\,\partial\sigma_d=\prod_{\tau_{d-1}<\sigma_d}\E\,\tau_{d-1}$, giving the expected values of the sufficient statistics as functions of the model parameters~$p_{\sigma_d}$.
  \item\label{com:surprise} The presence of the constraints on the existence of the boundaries of $d$-simplexes, their $d$-shells, \emph{in addition} to the constraints on the existence of $d$-simplexes themselves, may appear surprising at first. These stem from the conditional nature of the definition of these complexes. The models with these additional constraints removed, e.g.\ the canonical model over~$\cC_n$ with $\bx=\set{M_{s,d}}_{d=1}^{n-1}$, are still entirely unknown. This is not surprising, since such models appear to be combinatorially intractable~\cite{zuev_exponential_2015}.
  \item\label{com:hscm} Note that the constraints under which entropy is maximized in a hypercanonical model are generally intractable. The proofs that these constraints are the degree distributions in the dense and sparse HSCM graphs appear in \cite{chatterjee_random_2011} and \cite{hoorn_sparse_2018}, respectively. That is, the HSCM is the unbiased maximum-entropy model of random graphs with a given degree \emph{distribution}, versus an (expected) degree \emph{sequence} in the (S)CM.
  \item\label{com:hscm_d} The $d$-dimensional ((H)S)CM models have not been considered before. The proof that $\SCM_d(n,\bk)$ maximizes entropy subject to the expected $(d-1,d)$-degree sequence constraints is a straightforward adjustment of notations in the corresponding proof for a general canonical model~\cite{jaynes_information_1957}, so that it is omitted here for brevity.
  \item\label{com:z-hscm_d} The hypercanonical version of the model (Z-$\HSCM_d(n,\rho)$) has not been considered before. Efficient algorithms to sample from a generalized version of Z-$\CM_d(n,\bk)$ are considered in~\cite{young_construction_2017}.
\end{enumerate}

\section{Phase transitions}\label{sec:PT}
Phase transitions are very interesting and important phenomena occurring in random structures in statistical physics.
In this section we briefly recall the most fundamental types of phase transitions that have been studied in random graphs, and then describe their higher dimensional analogues in simplicial complexes.

\subsection{Homological Connectivity}
The very first theorem proved for random graphs in \cite{erdos_random_1959} was  the phase transition when the $G(n,M)$ random graph becomes connected. The following is the equivalent result for the $G(n,p)$ random graph. 

\begin{thm}
\label{thm:conn_gnp}
Let $G\sim G(n,p)$. Then
\[
    \limninf\prob{G\text{ is connected}} = \begin{cases}1 & np = \log n + w(n),\\
    0 & np = \log n -w(n),\end{cases}
\]
where $w(n)$ is any function satisfying $w(n)\to\infty$. In addition, if $np = \log n +c$ for $c\in (0,\infty)$, and denoting by $N_{\mathrm{comp}}$  the number of connected components in $G$, then
\[
(N_{\mathrm{comp}} -1) \xrightarrow{D} \pois{e^{-c}},
\]  
where $\xrightarrow{D}$ denotes convergence in distribution. This convergence implies that 
\[
    \limninf\prob{G \text{ is connected}} = e^{-e^{-c}}. 
\]
\end{thm}

The main idea behind the proof of Theorem \ref{thm:conn_gnp} is to show that around $np=\log n$, the random graph consists only of a giant connected component and isolated vertices. In that case, the phase transition for connectivity can be achieved by analysing the number of isolated vertices.

Note that in the language of homology, the event `$G$ is connected' can be phrased as `the $0$th homology group $H_0(G)$ is trivial', and that $N_{\mathrm{comp}}$ is equal to the $0$th Betti number $\beta_0(G)$.
Thus, it is tempting to try to generalize this phase transition for higher degrees of homology, and search for the point (value of $p$ or $np$) where the $k$th homology group $H_k$ becomes trivial.

\subsubsection{The random $d$-complex}

We start with the Linial-Meshulam random $d$-complex $Y_d(n,p)$. Recall that $G(n,p) = Y_1(n,p)$, and that Theorem \ref{thm:conn_gnp} studies $H_0$. Similarly, the following results studies $H_{d-1}$ in $Y_d(n,p)$.

\begin{thm}
\label{thm:conn_yd}
Let $Y\sim Y_d(n,p)$. Then
\[
    \limninf\prob{H_{d-1}(Y) = 0} = \begin{cases} 1 & np = d\log n + w(n),\\ 0 & np = d\log n - w(n),\end{cases}
\]
where $w(n)\to\infty$. 

In addition, if $np = d\log n + c$, $c\in (0,\infty)$, then
\[
    \beta_{d-1}(Y) \xrightarrow{D} \pois{\frac{e^{-c}}{d!}},
\]
which implies that
\[
    \limninf\prob{H_{d-1}(Y) = 0} = e^{-\frac{e^{-c}}{d!}}.
\]
\end{thm}
Note that taking $d=1$ in Theorem~\ref{thm:conn_yd}, nicely recovers the graph result in Theorem~\ref{thm:conn_gnp}. Thus, Linial and Meshulam named this phase transition `homological connectivity'.
The phase transition itself was proved first for $d=2$ and for homology in $\Z_2$ coefficients in \cite{linial_homological_2006}, and then for any $d$ and homology in $\Z_m$ coefficients  in \cite{meshulam_homological_2009}. The Poisson limit was proved a decade later \cite{kahle_inside_2016}.

Aside from the analogy between the statements of Theorems \ref{thm:conn_yd} and \ref{thm:conn_gnp}, the general idea behind the proof also shares some similarity. To prove Theorem \ref{thm:conn_yd} one can show that around $np = d \log n$, the only possible cycles in $H_{d-1}$\footnote{More accurately, the proof actually looks at cocycles in $H^{d-1}$.} are those generated by isolated $(d-1)$-simplexes. An `isolated' $(d-1)$-simplex is such that it is not included in any $d$-dimensional simplex. The fact that isolated $(d-1)$-simplexes yield nontrivial $H_{d-1}$ is relatively easy to prove. The much harder part here is to show that these are the \emph{only} possible cycles.

\subsubsection{The random flag complex}
\label{sec:conn_flag}

The random flag complex $X(n,p)$ differs from the random $d$-complex $Y_d(n,p)$ in two important ways. Firstly, $X(n,p)$ has random homology in all possible degrees $k$, rather than just $k=d$ and $k=d-1$ as in $Y_d(n,p)$. Secondly, note that in $Y_d(n,p)$ both $H_{d-1}$ and $H_d$ are \emph{monotone} -- $H_{d-1}$ is decreasing in $p$, while $H_d$ is increasing.
This is not the case for the flag complex, where except for $H_0$, none of the homology groups is monotone. The following result was proved by Kahle \cite{kahle_sharp_2014}.
\begin{thm}\label{thm:conn_xnp}
Let $X\sim X(n,p)$. Then,
\[
    \limninf\prob{H_k(X) = 0} = \begin{cases}1 & np^{k+1} = \param{\frac{k}{2}+1+\eps}\log n,\\
    0 & np^{k+1} = \param{\frac{k}{2}+1-\eps}\log n. \end{cases}
\]
\end{thm}
Note that here as well, taking $k=0$ agrees with the phase transition in Theorem \ref{thm:conn_gnp}. The phase transition here is also a consequence of the vanishing of the isolated $k$-simplexes. However, as oppose to the proofs in \cite{linial_homological_2006,meshulam_homological_2009} which  mainly consist of combinatorial arguments, the proof in \cite{kahle_sharp_2014} goes in a different way, employing \emph{Garland's method}. Briefly, Garland's method \cite{garland_p-adic_1973} is a powerful tool that allows ``breaking" the computation of homology into local pieces (the \emph{links} of the faces), and to invoke spectral graph arguments.

\subsubsection{The multi-parameter random complex}

Phase transitions are usually described as a rapid change of system behavior, in response to a small change in a system parameter value (e.g.~$p$ or $r$ above). Since the multi-parameter complex $X(n,\bp)$ has an infinite sequence of parameters, we cannot describe a suitable phase transition {\it per se}. Nevertheless, interesting results were presented in \cite{costa_large_2017},
showing that there are different regions within this high-dimensional parameter space, where in each region there is a single dominant homology group $H_k$.  

To describe the result we need a few definitions. Let $d>0$ and let $\alpha = (\alpha_1,\ldots,\alpha_d)\in \R^d_+$. Define
\[
\psi_k(\alpha) = \sum_{i=1}^d \binom{k}{i} \alpha_i,\quad k=1,\ldots,d.
\]
These functions are then used to define different domains
\[
\cD_k = \set{\alpha\in \R^d_+ : \psi_{k}(\alpha) < 1 < \psi_{k+1}(\alpha)},
\]
for $k=0,\ldots,d-1$, and where $\psi_0\equiv 0$. We also define $\cD_d = \set{\alpha : \psi_d(\alpha) < 1}$.
Finally, define
\[
    \tau_k(\alpha) = \sum_{i=1}^k(1-\psi_i(\alpha)),\qquad e(\alpha) = \min_{0\le i \le d}(1-\psi_i(\alpha)).
\]

\begin{thm}[\cite{costa_large_2017}]
Let $\bp = (p_1,p_2,\ldots,p_d,0,0,\ldots)$, be such that $p_i = n^{-\alpha_i}$, where $\alpha = (\alpha_1,\ldots,\alpha_d)\in \R^d$ may be a function of $n$, and $\limninf \alpha(n) = \alpha_*$.
Let $X\sim X(n,\bp)$, and suppose that $\alpha_* \in \cD_k$. Then, with high probability
\[
\beta_k(X) \approx n^{\tau_k(\alpha)}{(k+1)!},
\]
and for all $j\ne k$ we have
\[
\beta_j(X)  = O(n^{-e(\alpha)}\beta_k(X)).
\]
\end{thm}
In other words, in each $\cD_k$ the $k$-th Betti numbers is quite large, and all other Betti numbers are negligible compared to it. One can think of this results as a ``multi-parameter phase transition" so that the system behavior changes as one moves from one region $\cD_k$ to another. The proofs in \cite{costa_large_2016} relies mainly on counting faces and using Morse inequalities.

\subsubsection{Random geometric complexes}
\label{sec:conn_geom}

With respect to connectivity, the random geometric graph $G(\cP_n,r)$ behaves quite similarly to the $G(\cP_n,p)$ random graph. Denoting by $\Lambda = \omega_d nr^d$ the expected degree equal to the expected number of points in a ball of radius~$r$, with $\omega_d$ standing for the volume of the unit-ball in $\R^d$, the following was shown in \cite{penrose_longest_1997}.

\begin{thm}
\label{thm:conn_gnr}
Let $G\sim G(\cP_n,r)$. Then
\[
    \limninf\prob{G\text{ is connected}} = \begin{cases}1 & \Lambda = \log n + w(n),\\
    0 & \Lambda = \log n -w(n),\end{cases}
\]
where $w(n)$ is any function satisfying $w(n)\to\infty$.
\end{thm}
In other words, in both models connectivity is achieved once the expected degree is larger than $\log n$.

Moving to higher dimensions, the geometric models start to exhibit different behavior than the combinatorial ones. The main difference is the following. In $Y_d(n,p)$ and $X(n,p)$, when $p\to 1$, homological connectivity  describes the stage where homology becomes trivial. This is due to the fact that there is no structure underlying the complex. In the geometric complexes, the vertices are sampled over a metric space $S$, which might has its own intrinsic homology. Thus, for $r$ large enough, one should expect the homology of the random complex ($C(\cP_n,r)$ or $R(\cP_n,r)$) to ``converge" to the homology of $S$, rather than to vanish.
To account for this, the event we will refer to as the $k$-th homological connectivity here, is defined in \cite{bobrowski_homological_2019} as
\[
\cH_{k,r} := \set{\forall s\ge r: H_k(C(\cP_n,s))\cong H_k(S)}
\]
for the \cech complex, and similarly for the Rips complex, with `$C$' replaced by `$R$'. Note that $\cH_{k,r}$ is a monotone event (i.e.~$\cH_{k,r_1}\subset \cH_{k,r_2}$ for all $r_1\le r_2$). 

In the \cech complex, the phase transition in $\cH_{k,r}$ is by now fully understood. The following was proved in \cite{bobrowski_homological_2019}.

\begin{thm}\label{thm:conn_cnr}
Let $S=\T^d$, then for $1\le k \le d-2$, 
\[
\limninf \prob{\cH_{k,r}} = \begin{cases}
1 & \Lambda = 2^d(\log n + (k-1)\logg n + w(n)),\\
0 & \Lambda = 2^d(\log n + (k-1)\logg n - w(n)),
\end{cases}
\]
for any $w(n)\to\infty$.
In addition,
\[
\limninf \prob{\cH_{d-1,r}} =\limninf \prob{\cH_{d,r}} = \begin{cases}
1 & \Lambda = 2^d(\log n + (d-1)\logg n + w(n)),\\
0 & \Lambda = 2^d(\log n + (d-1)\logg n - w(n)).
\end{cases}
\]
\end{thm}
 While the proof there considers only the flat torus case, the results in \cite{bobrowski_random_2019} indicate that similar results could be achieved for any compact Riemannian manifold.
The key idea behind the proof of Theorem \ref{thm:conn_cnr} was to consider the evolution of complex ${C(\cP_n,r)}$ as $r$ is increased, and to look for \emph{critical faces}. By that we mean simplexes that facilitate changes in $H_k(C(\cP_n,r))$ when they first enter the complex. The analysis of critical faces employs Morse theory \cite{milnor_morse_1963,gershkovich_morse_1997}. The proof then consists of the following arguments:
\begin{enumerate}
    \item If $\Lambda$ is large enough ($\Lambda \gg \log n$), we have coverage, i.e.~$B(\cP_n,r) = \T^d$. Then, by the Nerve Lemma \cite{borsuk_imbedding_1948} we have that $H_k(C(\cP_n,r)) \cong H_k(\T^d)$.
    
    \item If $\Lambda  = 2^d(\log n+(k-1)\logg n + w(n))$ we can show that w.h.p.~no more critical faces that modify $H_k$ enter the complex. In other words, $H_k$ has reached its limit. From the previous argument the limit is $H_k(\T^d)$, and thus we conclude that $\cH_{k,r}$ holds. 
    
    \item If $\Lambda  = 2^d(\log n+(k-1)\logg n - w(n))$, critical faces still enter the complex and make chances in $H_k$. Thus $\cH_{k,r}$ does not hold (w.h.p.). 
\end{enumerate}

Similarly to Theorems \ref{thm:conn_gnp} and \ref{thm:conn_yd}, the results in \cite{bobrowski_homological_2019} also include a Poisson limit for the counting of the critical faces, when $\Lambda = 2^d(\log n + (k-1)\logg n + c)$. This result is the analogue of the isolated vertices in the random graph, since: (a) Critical faces are the `obstructions' to homological connectivity, in the same way that isolated vertices are for connectivity. (b) It can be shown that critical faces are indeed isolated when they first enter the complex. Since the exact definition of the critical faces require more details, we will not present the exact theorem here, but refer the reader to \cite{bobrowski_homological_2019,bobrowski2021poisson}.

Homological connectivity for the Rips complex $R(\cP_n,r)$ is still an open problem to date. The most recent result in this direction appeared in \cite{kahle_random_2011}. By using Discrete Morse Theory \cite{forman_users_2002}, Kahle was able to prove the following.
\begin{thm}
Let $S$ be a compact and convex subset of $\R^d$. Then
\[
    \mean{\beta_k(R(\cP_n,r))} = O\param{n\Lambda^k e^{-c\Lambda}}.
\]
In particular, if $\Lambda = \frac{1+\eps}{c}\log n$ then $\beta_k=0$ (w.h.p.).
\end{thm}
Note that this result is for compact and convex sets $S$, for which $H_k(S)=0$ for all $k\ge 1$. While this result is a significant step, it is still incomplete since: (a) It does not provide a sharp phase transition. In particular, there is no reason to believe that the constant $c$ provided by the proof, is the optimal one. (b) It does not apply to general manifolds.

\subsection{Emergence of homology}\label{sec:hom_emerge}

Another fundamental result proved from random graph concerns the appearance of cycles. The following was proved in \cite{erdHos1960evolution} (for the $G(n,M)$ model, but the statement for $G(n,p)$ is equivalent).

\begin{thm}\label{thm:cyc_gnp}
Let $G\sim G(n,p)$, then
\[
\limninf\prob{G \text{ contains  cycles}} = \begin{cases}
1 & np = c \ge 1,\\
\gamma(c) & np =c\in (0,1),\\
0 & np = o(1),
\end{cases}
\]
where $\gamma(c) = 1-\sqrt{1-c}\exp\param{\frac{c}{2}+\frac{c^2}{4}}$.
\end{thm}
Note that the event `$G$ contains cycles' can be phrased as $H_1(G) \ne 0$. In other words, Theorem \ref{thm:cyc_gnp} describes a phase transition for the emergence of $H_1$ in $G(n,p) = X_1(n,p)$.

In the following, we will review the results known to date about the emergence of $k$-cycles for the various random simplicial complex models.

\subsubsection{The random $d$-complex}

For the random $d$-complex, the following is an aggregate of a collection of works.
\begin{thm}
Let $Y\sim Y_d(n,p)$. There exists $c_d^*>0$ such that
\[
    \limninf\prob{H_{d}(Y) \ne 0} = \begin{cases} 
    1 & np = c > c^*_d,\\
    \gamma_d(c) & np = c \in (0,c^*_d),\\
    0 & np = o(1),
    \end{cases}
\]
where $\gamma_d(c) = 1-\exp\param{-\frac{c^{d+2}}{(d+2)!}}$. In addition,  the only $d$-cycles in $c\in (0,c_d^*)$ are generated by empty $(d+1)$-shells.
\end{thm}
 The first case was proved in \cite{aronshtam_collapsibility_2013}, the middle case was proved in \cite{linial_phase_2016}, and the last case was proved in \cite{kozlov_threshold_2010}. An equation defining the critical values can be found in \cite{aronshtam_collapsibility_2013}, while numerical approximations \cite{linial_phase_2016} yield $c_2^* \approx 2.754$, $c_3^* \approx 3.907, c_4^* \approx 4.962,\ldots, c_{1000}^* \approx 1001$. So roughly, it looks like $c_d^*\approx d+1$.

Another closely related phase transition is for collapsibility. Briefly, a $d$-complex is collapsible, if we can iteratively erase pairs of simplexes in dimension $d$ and $d-1$, without changing the topology of the complex. In \cite{aronshtam_threshold_2015} a phase transition for collapsibility was shown at critical values $c_d^{\mathrm{col}} < c_d^*$. 

\subsubsection{The random flag complex}
Recall that the homology $H_k(X(n,p))$ is non-monotone (in $p$). In Section \ref{sec:conn_flag} we saw that the largest value of $p$ such that $H_k(X(n,p))$ is non-trivial is when $np^{k+1} = \Theta(\log n)$. In this section we  look for the smallest possible value of $p$.

\begin{thm}
Let $X\sim X(n,p)$, and suppose that $\limninf np^k < \infty$.  Then,
\[
\limninf\prob{H_k(X) \ne 0} = \begin{cases} 1 & np^k \ge k+1+\eps,\\
0 & np^k \le n^{-\eps},
\end{cases}
\]
for any $\eps>0$.\end{thm}
The upper case was proved in \cite{kahle_sharp_2014}, while the lower case was proved in \cite{kahle_topology_2009}.
While this result is not as sharp as the previous ones, it still implies a phase transition when $np^k = \Theta(1)$. Together with Theorem \ref{thm:conn_xnp}, we have that $H_k(X)\ne 0$  for $p$ in the interval
\[
\sbrk{\param{\frac{k+1+\eps}{n}}^{\frac{1}{k}} , \param{\frac{\param{\frac{k}{2}+1-\eps}\log n}{n}}^{\frac{1}{k+1}}}.
\]

\subsubsection{Random geometric complexes}
The behavior of random geometric complexes is similar in many ways to that of the random flag complex. In particular, for each $k$ we have two phase transitions -- one for the emergence of the first $k$-cycles, and another for homological connectivity (Section \ref{sec:conn_geom}).
Here we present the former phase transition.

For the \cech complex, combining the results from \cite{kahle_random_2011,kahle_limit_2013}, we have the following.
\begin{thm}\label{thm:cyc_cnr}
For $1\le k \le d-1$,
\[
\limninf\prob{H_k(C_r(\cP_n)) \ne 0} = \begin{cases}
1 & n\Lambda^{k+1} = \omega(1),\\
\gamma_k(c) & n\Lambda^{k+1} = c\in (0,\infty),\\
0 & n\Lambda^{k+1} = o(1),
\end{cases}
\]
where $\gamma_k(c) = 1-\exp(-cA_k)$, for some constant $A_k>0$.
\end{thm}
Note that the transition occurs when $\Lambda = \Theta(n^{-\frac{1}{k+1}}) \to 0$. In other words, the first cycles appear in a sparse regime where the expected degree is small. Consequently, it was shown in \cite{kahle_random_2011} that in this regime the only $k$-cycles that appear are the smallest possible ones -- i.e., $(k+1)$-shells, which consist of $(k+2)$ vertices. The phase transition in Theorem \ref{thm:cyc_cnr} then essentially describes the appearance of these empty $(k+1)$-shells.

A similar result was proved for the Rips complex.
\begin{thm}\label{thm:cyc_rnr}
For any $k\ge 1$, 
\[
\limninf\prob{H_k(R_r(\cP_n)) \ne 0} = \begin{cases}
1 & n\Lambda^{2k+1} = \omega(1),\\
\tilde\gamma_k(c) & n\Lambda^{2k+1} = c\in (0,\infty),\\
0 & n\Lambda^{2k+1} = o(1),
\end{cases}
\]
where $\tilde\gamma_k(c) = 1-\exp(-c\tilde A_k)$, for some constant $\tilde A_k>0$.
\end{thm}
There are two main differences between Theorems \ref{thm:cyc_cnr} and  \ref{thm:cyc_rnr}. Firstly, due to the Nerve Lemma \cite{borsuk_imbedding_1948}, the \cech complex on $\T^d$ (as well as any subset of $\R^d$) cannot have any $k$-cycles for $k>d$. The Rips complex, on the other hand, may have $k$-cycles for any $k$.
Secondly,  the transition here occurs when $\Lambda = \Theta\param{n^{-\frac{1}{2k+1}}}$, as opposed to $\Lambda = \Theta\param{n^{-\frac{1}{k+1}}}$ in Theorem \ref{thm:cyc_cnr}. This difference stems from the fact that empty shells cannot appear in the Rips complex. Thus, it can be shown that the smallest possible $k$-cycle is achieved by generating an empty cross-polytope, which consists of $2k+2$ vertices. 

\subsection{Percolation-related phenomena}

Percolation theory studies the emergence of infinite or ``giant" connected component in various random media. The study of higher-dimensional analogues is at a preliminary stage, and we describe here the most recent results.

In the case of random graphs, giant components are commonly defined in terms of the number of vertices. Specifically, for both the $G(n,p)$ and $G(\cP_n,r)$ models, by a ``giant component" we refer to a component that consists of $\Theta(n)$ many vertices. 
As opposed to connectivity which we discussed earlier, this notion does not have a simple analogue in higher dimensions. In this section we present two different notions -- one of the random $d$-complex and another for the random \cech complex. 

\subsubsection{Random $d$-complex}

Starting with the $G(n,p)$ graph, the following was proved in \cite{erdHos1960evolution}. \begin{thm}\label{thm:gnp_giant}
Let $G\sim G(n,p)$, and $np = c\in(0,\infty)$. Denote by $L_n$ the largest connected component in $G$. Then with high probability
\[
L_n = \begin{cases} \Theta(n) & c > 1,\\
\Theta(n^{2/3}) & c=1,\\
\Theta(\log n) & c< 1.
\end{cases}
\]
\end{thm}
In other words, a giant component appears for $c>1$. Further, it can be shown that in this case, the giant component is unique.

In order to generalize the emergence of the giant component to higher dimension, the notion of \emph{shadow} was introduced in~\cite{linial2019extremal}.
Suppose that $X$ is a $d$-dimensional simplicial complex with a complete $(d-1)$-skeleton. The shadow of $X$, denoted $\mathrm{SH}(X)$, is the set of all $d$-dimensional faces $\sigma$ such that (a) $\sigma\not\in X$, and (b) adding $\sigma$ to $X$ generates a new $d$-cycle.

In the case of the graph (i.e., $d=1$) the shadow are all the edges not in the graph, such that adding them generates a new cycle. In that case, it is easy to show that when $c>1$ we have $|\mathrm{SH}(X)| = \Theta\param{\binom{n}{2}}$, so that the shadow has a \emph{positive density}, and we say that the graph has a \emph{giant shadow}. In \cite{linial_phase_2016}, this phenomenon was generalized for higher dimensions.

\begin{thm}
Let $Y\sim Y_d(n,p)$, and $np = c\in (0,\infty)$. Then with high probability,
\[
|\mathrm{SH(X)}| = \begin{cases} \Theta\big(\binom{n}{d+1}\big) & c > c_d^*,\\
\Theta(n) & c < c_d^*,\end{cases}
\]
where $c_d^*$ is the critical value defined in Section \ref{sec:hom_emerge}.
\end{thm}
In other words, when $c>c_d^*$ a giant shadow emerges in the random $d$-complex, similarly to what happens in the graph case.

\subsubsection{Random geometric complexes}\label{sec:giant_cycles}

Similarly to the $G(n,p)$ model, the random geometric graph $G(\cP_n,r)$ exhibits a phase transition for the emergence of a giant component. This phase transition occurs when the expected degree $\Lambda$ is constant, and the corresponding critical value is denoted as $\lambda_c$.

In the context of geometric complexes, a different generalization for percolation was studied recently \cite{bobrowski_homological_2020,bobrowski2020homological}. Consider the random \cech complex $C(\cP_n,r)$ defined in Section \ref{sec:models}, where $\cP_n$ is a homogeneous Poisson process on the $d$-torus $\T^d$. One of the nice facts about the torus, is that it has non-trivial homology for all $0\le k\le d$, and in particular $\beta_k(\T^d) = \binom{d}{k}$. In \cite{bobrowski_homological_2020} the notion of a ``giant $k$-cycle" was introduced, by which we mean any $k$-cycle in $C(\cP_n,r)$ that corresponds to any of the nontrivial cycles in $H_k(\T^d)$.

To define this more rigorously, we consider the union of balls $B(\cP_n,r)$, and the fact that $H_k(B(\cP_n,r))\cong H_k(C(\cP_n,r))$. We then take the inclusion map $i:B(\cP_n,r)\hookrightarrow \T^d$, and its induced map $i_*:H_k(B(\cP_n,r))\to H_k(\T^d)$. By a \emph{giant $k$-cycle} we refer to all nontrivial elements in the image $\im(i_*)$. Finally, we define the events
\[
    {\cal E}_k := \set{\im(i_*) \ne 0}, \quad {\cal A}_k = \set{\im(i_*) = H_k(\T^d)}.
\]
In other words, $\cE_k$ is the event that some giant cycles exist, while $\cA_k$ is the event that all of them are present in $B(\cP_n,r)$ (and correspondingly in $C(\cP_n,r)$). The following was proved in \cite{bobrowski_homological_2020}.

\begin{thm}\label{thm:hom_perc}
Suppose that $\Lambda = \omega_d nr^d = \lambda \in (0,\infty)$. 
\begin{enumerate}
    \item There exist $\lambda_1^0 \le \lambda_2^0\ldots\le\lambda_{d-1}^0$, such that if $\lambda < \lambda_k^0$ then,
\[
    \prob{{\cal A}_k} \le \prob{\cE_k} \le e^{-C_k^0 n^{1/d}},
\]
for some $C_k^0 > 0$.
\item There exist $\lambda_1^1 \le \lambda_2^1\ldots\le\lambda_{d-1}^1$, such that if $\lambda > \lambda_k^1$ then,
\[
    \prob{\cE_k} \ge \prob{\cA_k} \ge 1-e^{-C_k^1 n^{1/d}},
\]
for some $C_k^1 > 0$.
\end{enumerate}
In addition, we have $\lambda_k^0 \le \lambda_k^1$ for all $k$, and $\lambda_1^0 = \lambda_1^1 = \lambda_c$.
\end{thm}
In other words, Theorem \ref{thm:hom_perc} suggests that there is a phase transition describing the emergence of the giant $k$-cycles. 
However in order to complete the proof, one needs to show that $\lambda_k^0 = \lambda_k^1$ for all $k$. 

At this stage the giant shadow and the giant cycle phenomena are incomparable. However, it should be noted that both occur in the regime where the expected degree ($np$ or $\Lambda$) is finite, similarly to the giant component. It is an interesting question whether these are merely two ways to view the same phenomenon, or they describe completely different structures. Finally, we should note that other ideas for extending percolation-type phenomena to higher dimension exist in the literature.

\subsection{The fundamental group}
Given a space $X$, the fundamental group $\pi_1(X)$ represents the space of equivalence classes of loops in $X$, where the equivalence is based on homotopy (a continuous transformation of one loop into the other). The space $X$ is simply connected if and only if $\pi_1(X)$ is trivial. 

Note that the first homology group  $H_1(X)$  also represents loops, but under a different equivalence relation (being a boundary). In fact, it can be shown that $H_1(X)$ is an abelianization of $\pi_1(X)$. Thus, it is possible that $H_1(X)$ is trivial while $\pi_1(X)$ is not. Consequently, the vanishing threshold for $\pi_1$ is larger than the threshold for $H_1$ as the following statements show.

We start with the random $d$-complex. Here the only relevant dimension is $d=2$, since otherwise the fundamental group is trivial. The following was proved in \cite{babson_fundamental_2011,luria2018simple}.
\begin{thm}
Let $Y\sim Y_2(n,p)$. Then,
\[
\limninf \prob{\pi_1(Y) = 0} = \begin{cases}1 &  np = c n^{\frac{1}{2}}, \\ 0 & np = n^{\frac{1}{2}-\eps},\end{cases}
\]
where $c$ is any constant bigger than $\sqrt{3^3/4^4}$.  
\end{thm}
Further, when $np = n^{1/2-\eps}$ it is shown in \cite{babson_fundamental_2011} that the fundamental group is hyperbolic. 
It is also conjectured in \cite{luria2018simple} that in fact $c = \sqrt{3^3/4^4}$ is the exact sharp threshold for simple connectivity.

For the random flag complex, the existing results \cite{babson_fundamental_2012}
 are even less sharp.
\begin{thm}
Let $X\sim X(n,p)$, then for any 
\[
    n^{\frac{1}{2}+\eps} \le np \le n^{\frac{2}{3}-\eps},
\]
with high probability $\pi_1(X)$ is nontrivial and hyperbolic.
\end{thm}

Note that the algebraic structure of the fundamental group is more intricate than $H_1(X)$, and thus several other properties have been studied in the literature. For example: its geometric and cohomological dimension \cite{costa_fundamental_2015}, torsion \cite{costa_geometry_2015}, property T \cite{kahle_sharp_2014}, freeness \cite{newman2020freeness}, and more.

\section{Limit theorems}\label{sec:LT}
The results in the previous section focused on a rather qualitative analysis of phase transitions.
Another direction of study is to consider some numerical marginals of homology (e.g., the Betti numbers) and find their limiting distribution.
In this section we will briefly review some of the results in this direction.

We will not discuss the proofs here, however we note that in most cases they rely on \emph{Stein's method} \cite{stein1986approximate,ross2011fundamentals}. Briefly, this method allows one to prove central limit theorems and convergence to the normal and Poisson distributions in cases that involve sums of random variables that are \emph{not} independent.

\paragraph{Notation.}
In this section we will use $a_n\approx b_n$ to denote that $a_n/b_n \to 1$ as $n\to\infty$, and $a_n\ll b_n$ to denote that $a_n/b_n \to 0$. In addition, $\xrightarrow{D}$ stands for convergence in distribution, and $N(\mu,\sigma^2)$ for the normal distribution with mean $\mu$ and variance $\sigma^2$.

\subsection{Betti numbers}
As we have seen in Section \ref{sec:hom_emerge}, the $k$-th homology of $X(n,p)$ is nontrivial roughly when $p$ is between $n^{-1/k}$ and $n^{-1/(k+1)}    $ (up to constants). It was shown in \cite{kahle_topology_2009}
that when $n^{-1/k} \ll p \ll n^{-1/(k+1)}$ then
\[
\mean{\beta_k(X(n,p))} \approx\frac{n^k}{(k+1)!} p^{\binom{k+1}{2}}.
\]
In addition, the following central limit theorem was proved in \cite{kahle_limit_2013}.
\begin{thm}
Let $X\in X(n,p)$, and let $n^{-1/k} \ll p \ll n^{-1/(k+1)}$. Then
\[
    \frac{\beta_k(X)-\mean{\beta_k(X)}}{\sqrt{\var{\beta_k(X)}}} \xrightarrow{D} N(0,1).
\]
\end{thm}

In the geometric complexes case, the results usually split according to the expected degree $\Lambda$. In the case where $\Lambda \to 0$ (sometimes called the sparse regime), the homology is dominated by empty shells, and thus it is relatively easy to compute the Betti numbers. For example, it was proved in \cite{kahle_random_2011} that 
\[
    \mean{\beta_k(C(\cP_n,r))} \approx A_k n \Lambda^{k+1},
\]
where $A_k>0$ is a known constant~\cite{kahle_random_2011}.
Note that the limit of the right hand side can be zero, finite, or infinite. Consequently we have three possible limits. The following was proved in \cite{kahle_limit_2013}.
\begin{thm}
Let $C\sim C(\cP_n,r)$, and suppose that $\Lambda\to 0$.
\begin{enumerate}
    \item If $n\Lambda^{k+1}\to 0$ then $\beta_k(C) = 0$ with high probability.
    \item If $n\Lambda^{k+1} = \alpha\in (0,\infty )$ then $\beta_k(C)\xrightarrow{D} \pois{A_k\alpha}$.
    \item If $n\Lambda^{k+1}\to\infty$ then $\frac{\beta_k(C) - \mean{\beta_k(C)}}{\sqrt{\var{\beta_k(C)}}} \xrightarrow{D} N(0,1)$.
\end{enumerate}
\end{thm}
Similar results hold for the Rips complex as well, just with $n\Lambda^{2k+1}$ replacing $n\Lambda^{k+1}$ \cite{kahle_limit_2013}.

The regime where $\Lambda$ is finite and non-vanishing is sometimes referred to as the \emph{thermodynamic limit}. It can be shown that most of the cycles in either the \cech or Rips complexes are generated in this regime. However, exact counting here is much more difficult. It was shown in \cite{kahle_random_2011} that $\mean{\beta_k(C(\cP_n,r))} = \Theta(n)$ (and the same for the Rips). However the exact expectation is not known. 

Nevertheless, one can prove laws of large numbers as well as central limit theorems for the Betti number functionals.
This type of results was first introduced in \cite{yogeshwaran_random_2016}, and have been improved over the years \cite{hiraoka_limit_2018,krebs2019asymptotic,trinh2019central}. The main tools used are stabilization methods for random point processes (e.g.~\cite{penrose2001central}). The result is as follows.

\begin{thm}\label{thm:betti_clt}
Let $C\sim C(\cP_n,r)$, and suppose that $\Lambda = \lambda\in (0,\infty)$. Then $$\frac{\beta_k(C) - \mean{\beta_k(C)}}{\sqrt{n}} \xrightarrow{D} N(0,\sigma^2(\lambda)),$$
for some $\sigma^2(\lambda)>0$.
\end{thm}
A similar result applies to the Rips complex as well. In addition in \cite{trinh2019central,goel2019strong} the binomial process $\cX_n$ was studied as well, and similar limiting theorems were proved.

In the dense regime, i.e., $\Lambda\to\infty$ there are very few results. In particular, even the scale of the expected Betti numbers is not known. The only case analyzed is the Poisson limit when $\Lambda \approx \log n$ described earlier, as part of the homological connectivity phenomenon.

\subsection{Topological types}
In \cite{auffinger_topologies_2020} a completely different type of limit theorem was proved. There, the goal was to study the distribution of the different topological types that may appear in a random \cech complex. Let $X$ be a simplicial complex, and let $\cC(X)$ be the set of all connected components of $X$ (so that $\beta_0(X) = |\cC(X)|$).
Define the empirical measure
\[
\mu_X := \frac{1}{\beta_0(X)}\sum_{C\in \cC} \delta_{[C]},
\]
where $\delta$ stands for the Dirac delta measure, and $[C]$ stands for the equivalence class of all components homotopy equivalent to $C$ (i.e.~so that one component can be ``continuously transformed" into the other, see \cite{hatcher_algebraic_2002}).

\begin{thm}
Let $C\sim C(\cP_n,r)$, and denote $\hat\mu_n :=\mu_{C}$. Then the random measure $\hat\mu_n$ converges in probability to a universal probability measure.
The support of the distribution is the set of all connected (finite) \cech complexes in $\R^d$.
\end{thm}
Note that universality here means that the limit is the same regardless of the underlying manifold $S$.

\subsection{Persistent homology}\label{sec:persistent}

 In the context of Topological Data Analysis (TDA), one is often more interested in studying \emph{persistent} homology rather than just the (static) homology.
In this section we consider the persistent homology over  filtrations of either \cech or Rips complexes, namely $\set{C(\cP_n,r)}_{r=0}^{\infty}$ or $\set{R(\cP_n,r)}_{r=0}^{\infty}$, respectively.

\subsubsection{Limit theorems}

One useful quantity that can be extracted from persistent homology are the \emph{persistent Betti numbers}. Briefly, for any $s\le t$, we denote by $\beta_k^{(s,t)}$ the number of cycles born before radius $s$ and die after radius $r$ (see \cite{hiraoka_limit_2018} for a formal definition). In \cite{hiraoka_limit_2018} the following central limit theorem was proved.

\begin{thm}\label{thm:pbetti_clt}
Let $C\sim C(\cP_n,r)$, and suppose that $ns^d = \alpha$ and $nt^d = \beta$, for some $\alpha \le \beta \in (0,\infty)$. Then
\[
\frac{\beta_k^{(s,t)}(C) - \mean{\beta_k^{(s,t)}(C)}}{\sqrt{n}} \xrightarrow{D} N(0,\sigma^2(\alpha,\beta)),
\]
for some $\sigma^2(\alpha,\beta)>0$.
\end{thm}
Note that this theorem describes the persistent Betti numbers in the thermodynamic limit, where most cycles are born and die. A similar theorem applies for $R(\cP_n,r)$, and in fact \cite{hiraoka_limit_2018} present this central limit theorem for an even larger general class of geometric complexes. In \cite{krebs2019asymptotic} Theorem \ref{thm:pbetti_clt} was further extended to a multidimensional central limit theorem (i.e.\ for a sequence $(\alpha_1,\beta_1),\ldots, (\alpha_m,\beta_m)$). Finally, note that Theorem \ref{thm:betti_clt} is in fact a special case of Theorem \ref{thm:pbetti_clt} in the case where $\alpha=\beta$.

Recall that a common representation for persistent homology is via a \emph{persistence diagram}, which is essentially a finite subset of $\Delta = \set{(x,y) : y\ge x} \subset \R^2$, where the coordinates of each point correspond to the birth and death times of a persistence interval.
Thus, persistence diagrams attached to random finite simplicial complexes, can be thought of as a random subset of $\Delta$, or equivalently a random discrete measure on $\Delta$. One can then employ the theory of random sets and random measures to study the limit of such diagrams. 

Consider the random \cech filtration $\set{C(\cP_n,r)}_{r=0}^\infty$, and let $\mu_k(n)$ be the discrete measure representing the persistence diagram in degree $k$. The main result in \cite{hiraoka_limit_2018} shows that there exists a limiting Radon measure on $\Delta$, 
\[
    \mu_k = \limninf \frac{\mu_k(n)}{n}.
\] The same is true for the Rips filtration as well. This statement is essentially a law of large numbers for the persistence diagram as a Radon measure in $\R^2$.

A different approach was taken in \cite{owada2020convergence}, where persistence diagrams were studied from the point of view of set-topology. In that case it was shown that each persistence diagram can be roughly split into three different areas. The limit of the persistence diagram in the \emph{sparse} region is empty. In the \emph{intermediate} region it has a limiting Poisson process distribution. Finally, in the dense region, the persistence diagram converges to a solid (nonrandom) 2d region.

\subsubsection{Maximal cycles}

Another interesting and significant type of limit for persistence diagrams was presented in \cite{bobrowski_maximally_2017}. There, it was assumed that $S$ is a unit box (though the results should not depend on that). The quantity considered was the maximal death/birth ratio among all persistent cycles of degree $k$, denoted $\Pi_k$. 

\begin{thm}\label{thm:max_cyc}
Let $\Pi_k$ be the maximally persistent cycle for either the Rips or \cech complex. Then there exist $A_k<B_k$ such that with high probability
\[
    A_k\param{\frac{\log n}{\logg n}}^{1/k} \le \Pi_k \le B_k\param{\frac{\log n}{\logg n}}^{1/k}.
\]
\end{thm}
Theorem \ref{thm:max_cyc} provides the scaling for the multiplicative persistence value in a random geometric complex. Since the unit box has no intrinsic homology, this result can be thought of as describing the homology of \emph{noisy} cycles, i.e., cycles that are not part of the underlying topology. It can be shown that the same result applies for other compact spaces $S$, as long as we do not count cycle that belong to $H_k(S)$. Recalling Section \ref{sec:giant_cycles}, such ``signal" cycles (those in $H_k(S)$) are born when $\Lambda$ is constant, implying that the birth time is $r\propto n^{-1/d}$. On the other hand, the signal cycles die when $r$ is a constant (depends on $S$ and not on $n$). Therefore, for the signal cycles we have $(\mathrm{death}/\mathrm{birth})\propto n^{1/d}$.  In other words, taking  death/birth as  a measure of persistence, Theorem \ref{thm:max_cyc} shows that asymptotically noisy cycles are significantly less persistent than the signal ones, and one should be able to differentiate between them, assuming a large enough sample.

\section{Other directions}

In Sections \ref{sec:PT} and \ref{sec:LT} we presented some of the most fundamental results proved in the mathematical literature on random simplicial complexes.
There are numerous important studies that we omit for space reasons, but would like to mention them briefly here.
\begin{itemize}
    \item {\bf Spectrum and expansion:} Similarly to graphs, one can define an adjacency and Laplacian operators for simplicial complexes, as well as different types of expansion properties. These were studied mainly in the the context of the random $d$-complex $Y_d(n,p)$~\cite{gundert_eigenvalues_2016,dotterrer_coboundary_2012,lubotzky_random_2019,knowles_eigenvalue_2017}.
    \item {\bf Functional limit theorems:}
    One can examine functionals such as the Betti numbers and the Euler characteristic in a dynamic setting, and seek a limit in the form of a stochastic (Gaussian) process.
    In \cite{owada2020convergence,thomas2021functional}, geometric complexes were studied, and the dynamics was the growing connectivity radius in the complex. The results show that the limiting process is indeed Gaussian.
    In \cite{thoppe2016evolution}, the random flag complex was studied. Here, the set of edges used to construct the complex turns on and off, following stationary Markov  dynamics. The limiting processes are shown to be  Ornstein-Uhlenbeck (Gaussian) processes.
    These results were further extended to the multiparameter complex \cite{owada2021limit}.
    \item {\bf Spanning acycles:} Similarly to the study of spanning trees in  graphs, one can consider \emph{spanning acycles} in simplicial complex. Suppose that $K $ is a complete $(d-1)$-skeleton on $n$ vertices. A set of $d$-faces $S$ is considered a spanning acycle if  $\beta_{d-1}(K\cup S) = \beta_d(K\cup S) = 0$.
    In other words, adding $S$ to $K$ kills all the existing $(d-1)$-cycles (hence ``spanning''), while not generating any new $d$-cycles (hence ``acycle''). This definition invites the study of spanning acycles in the random $d$-complex.
    In \cite{skraba_randomly_2017,hiraoka_minimum_2017} the weights of faces in the minimal spanning acycle and their connection to the persistence diagram were studied. 
    
\end{itemize}

\section{Future directions}\label{sec:conclusion}

We hope our review illustrates well that the recent progress in the fundamental study of random simplicial complexes in mathematics has been extremely rapid. As far as applications to real-world systems in network/data science are concerned, the progress has been as rapid, but the field is definitely less mature that its $1$-dimensional counterpart dealing with applications of random graph models. There, a collection of important and ubiquitous properties of real-world networks abstracted as graphs have been long identified and well studied. Such properties include sparsity, scale-free degree distributions, degree correlations, clustering, small-worldness, community structure, self-similarity, and so on~\cite{barabasi_network_2016,newman_networks_2018}. Consequently, there have been an impressive body of work on null models that reproduce these properties in the statistically unbiased manner, and on growing network models that attempt to shed light on possible mechanisms that might lead to the formation of these properties in real-world systems~\cite{barabasi_network_2016,newman_networks_2018}.

When the same systems are modeled as simplicial complexes or hypergraphs, the list of important properties one should be looking at, is much less understood. Consequently, the world of models dealing with such properties is much less explored, as it is not entirely clear what properties one should be most concerned with in the first place.

In addition to its relative infancy, the other reason why the world of random complex models has not yet been navigated as much, is that models of even simplest higher-dimensional properties tend to be much more complicated, not necessarily at the technical level as much as at the conceptual level.

Consider the null models of the degree distribution, i.e.\ the configuration models, for example. Since, as opposed to graphs, there is not one but ${n\choose2}$ notions of degree in a simplicial complex of size~$n$ ($d$-degrees of $d'$-simplexes for any pair of dimensions $d,d'$, $0\leq d'<d\leq n-1$), it may not be immediately clear what degrees to focus on. The degrees that the $\mathrm{(S)CM}_d(n,\bk)$ and Z-$\mathrm{(S)CM}_d(n,\bk)$ models reproduce in Sections~\ref{sec:scm} and~\ref{sec:zcm}, are the $d$-degrees of $(d-1)$- and $0$-dimensional simplexes, respectively, with very different constraints imposed on the underlying simplicial structure.

In that regard, these two models are the two simplest extremal versions of a great variety of configuration models we can think of. Such models can reproduce sequences of $d$-degrees of $d'$-simplexes for any pair of dimensions~$d,d'$, in which case the defining equations~(\ref{eq:scmd-p}, \ref{eq:scmd-k}) remain the same, except that simplexes $\tau$ in those equations are now understood to be of dimension~$d'$. We may also want to jointly reproduce different combinations of degree sequences of different dimensions~$d,d'$ satisfying structural constraints. Such models, based on generalized degree sequences observed in real data, could serve as more adequate and more realistic null models of the data. The first step in this direction has been made in~\cite{young_construction_2017}.

Another example is homology, which we discussed at length in this chapter. Homology is known to be a powerful tool to capture information about the qualitative structure of simplicial complexes, and its theory in the random setting is constantly growing. However, as a tool to analyze high-dimensional networks, it is not yet as clear what part of the information provided by homology is useful and for what purposes. For example, while $H_0$ (connected components) and $H_1$ (cycles) are quite intuitive, the role of $H_2$ as a network descriptor is not as clear. Neither is it entirely clear what useful information is contained in the number of $k$-cycles, their physical size, their persistent lifetime, etc. So while the theoretical homological analysis is fruitful and highly interesting, it tends to be a challenge to see how and when to use it in application to real-world networks.
 
The degree distribution and homology are certainly just two out of many properties that the future work may find interesting and important to model in more realistic models of random simplicial complexes. For now, the list of such models is quite rudimentary for the reasons above, so that it is not surprising that rigorous mathematical results are available only for very basic models. Yet we have seen that even for these simplest models, the spectrum of available results is extremely rich and complex.

We have also seen that these results mostly address direct problems, whereas in dealing with real-world data one often faces inverse problems. For example, the results presented in Section~\ref{sec:persistent} tell us facts about the persistence diagrams in ``laboratory-controlled'' settings, i.e.\ for given complexes in given spaces with known topology---a direct problem. An inverse problem would be to make any statistical inferences about the space topology from a persistence diagram coming from real-world data. Some progress in this direction has been made over the past decade \cite{niyogi2011topological,fasy2014confidence,chazal2017robust,bobrowski2017topological,reani2021cycle}, but this is still an active and important area of research to pursue.

Finally, we comment on one potentially very interesting direction of future research in mathematics of random complexes. One of the most fundamentally important recent achievements in mathematics of random graphs is the development and essential completion of the theory of limits of dense graphs known as graphons~\cite{lovasz_large_2012,janson_graphons_2013}. One of the main results in that theory is that, even though there are many very different notions of graph convergence, they all are equivalent for dense graphs, and if a family of dense graphs converge, they converge to a (random) graphon~\cite{diaconis_graph_2008}. As mentioned in Section~\ref{sec:general-random}, a graphon is an integrable connection probability function $W:[0,1]^2\to[0,1]$ modulo a certain equivalence relation, and a graphon-based random graph of size $n$---known as a $W$-random graph in the graphon theory---is $G\sim G(n,\hr)$, where $r_{ij}=W(X_i,X_j)$ and $X_i\sim U(0,1)$.
Can an analogous theory be constructed---or discovered, depending on one's philosophical view---for the limits of dense complexes, e.g.\ \emph{complexons}? Complexons may probably be related to \emph{hypergraphons}~\cite{gowers_hypergraph_2007,elek_measure_2012,zhao_hypergraph_2015,balasubramanian_nonparametric_2021}, the limits of dense hypergraphs, via the lower and upper complexes $\uZ(n,\hbp)$ and $\oZ(n,\hbp)$ discussed in Section~\ref{sec:lower-upper}.

\end{document}